\documentclass[12pt]{article}

\usepackage{amsmath}
\usepackage{amssymb}
\usepackage[curve,matrix,arrow,cmtip]{xy}
\NoComputerModernTips

\newtheorem{theorem}{\indent Theorem}[section]
\newtheorem{lemma}[theorem]{\indent Lemma}
\newtheorem{prop}[theorem]{\indent Proposition}
\newtheorem{propA}{\indent Proposition}
\newtheorem{cor}[theorem]{\indent Corollary}
\newtheorem{corA}[propA]{\indent Corollary}
\newtheorem{conj}[theorem]{\indent Conjecture}
\newtheorem{remark}[theorem]{\indent Remark}

\newtheorem{assumption}[theorem]{\indent Assumption}
\newtheorem{definition}[theorem]{\indent Definition}

\newcommand{\rhobar}{{\overline{\rho}}}

\advance\textheight by .1cm

\def\ublank{\underline{\phantom{n}}}
\DeclareMathOperator{\image}{im}
\DeclareMathOperator{\kernel}{ker}
\def\unr{{\mathrm{unr}}}
\def\BF{\mathbf{F}}
\def\BG{\mathbf{G}}
\def\BC{\mathbf{C}}
\def\BA{\mathbf{A}}
\def\BQ{\mathbf{Q}}
\def\BN{\mathbf{N}}
\def\BT{\mathbf{T}}
\def\BZ{\mathbf{Z}}
\def\eps{\varepsilon}
\def\ad{\mathrm{ad}^0(\rhobar)}
\def\Ad{\mathrm{ad}(\rhobar)}
\def\GL{\mathrm{GL}}
\def\ab{\mathrm{ab}}

\def\SL{\mathrm{SL}}

\def\Frob{\mathrm{Frob}}
\def\Gal{\mathrm{Gal}}

\def\Hom{\mathrm{Hom}}

\def\rank{\mathrm{rank}}

\DeclareMathOperator{\degr}{deg}
\DeclareMathOperator{\spec}{sp}
\DeclareMathOperator{\Ind}{Ind}

\DeclareMathOperator{\Diag}{Diag}
\DeclareMathOperator{\Sets}{\mathbf{Sets}}

\newcommand{\Ulambda}{\underline{\lambda}}
\newcommand{\Um}{{\underline{m}}}
\def\CA{\mathcal{A}}
\def\CH{\mathcal{H}}
\def\CC{\mathcal{C}}
\def\CD{\mathcal{D}}

\def\CO{\mathcal{O}}

\def\Fm{\mathfrak{M}}
\def\FM{\mathfrak{m}}

\newcommand{\bC}{\mathbf{C}}
\newcommand{\ObC}{\widetilde{\bC}}
\newcommand{\OpbC}{\widetilde{\mathbb{C}}\hspace*{.025em}}
\newcommand{\pbC}{\bC}

\def\theenumi{(\roman{enumi})}

\def\p@enumii{\theenumi}

\def\longto{\longrightarrow}
\def\into{\hookrightarrow}

\def\longinto{\lhook\joinrel\longrightarrow}

\newbox\mybox

\def\arrover#1{\mathrel{
       \setbox\mybox=\hbox spread 1.4em{\hfil$\scriptstyle#1$\hfil}
       \vbox{\offinterlineskip\copy\mybox
             \hbox to\wd\mybox{\rightarrowfill}}}}
\def\larrover#1{\mathrel{
       \setbox\mybox=\hbox spread 1.4em{\hfil$\scriptstyle#1$\hfil}
       \vbox{\offinterlineskip\copy\mybox
             \hbox to\wd\mybox{\leftarrowfill}}}}

\def\ontoover#1{\mathrel{
       \setbox\mybox=\hbox spread 1.4em{\hfil$\scriptstyle#1$\hfil}
       \vbox{\offinterlineskip\copy\mybox
             \hbox to\wd\mybox{\rightarrowfill\hskip-2.8mm
                               $\rightarrow$}}}}
\def\leftontoover#1{\mathrel{
       \setbox\mybox=\hbox spread 1.4em{\hfil$\scriptstyle#1$\hfil}
       \vbox{\offinterlineskip\copy\mybox
             \hbox to\wd\mybox{$\leftarrow$\hskip-2.8mm
                               \leftarrowfill}}}}
\def\onto{\ontoover{\ }}

\newif\ifnormalesBeweisEnde
\def\BeweisEnde{\EndOfBeweis
                \global\normalesBeweisEndefalse}
\newenvironment{Proofof}[1]%
  {\vskip 0ex plus 0.15ex minus 0.6ex \pagebreak[1]
   \global\normalesBeweisEndetrue
   \trivlist
   \item[\hskip\labelsep {\textsc{Proof} \rm of #1}:]}%
  {\ifnormalesBeweisEnde \EndOfBeweis \fi
   \endtrivlist
   \vskip 0ex plus .15ex minus 0.6ex \pagebreak[2]}
\newenvironment{Proof}%
  {\vskip 0ex plus 0.15ex minus 0.6ex \pagebreak[1]
   \global\normalesBeweisEndetrue
   \trivlist
   \item[\hskip\labelsep \textsc{Proof}:]}%
  {\ifnormalesBeweisEnde \EndOfBeweis \fi
   \endtrivlist
   \vskip 0ex plus .15ex minus 0.6ex \pagebreak[2]}
\def\EndOfBeweis{\hskip .5em \vrule width 1.0ex height 1.0ex depth 0.3ex}

\newcommand{\notdiv}{\mathopen{\mathchoice
             {\not{|}\,}
             {\not{|}\,}
             {\!\not{\:|}}
             {\not{|}}
             }}

\begin{document}
\hfuzz=2pt

\title{Mod $\ell$ representations of arithmetic fundamental groups II \\
{\small (A conjecture of A.J. de Jong)}}

\author{Gebhard B\"ockle \thanks{G.B would like to thank the
TIFR for its hospitality in the Summer of 2002 during which the first
decisive steps in this work were made and also the ETH for its
inspirational environment and for its generosity in providing a three
year post-doctoral position.} \\ 
Chandrashekhar Khare \thanks{Some of the work on this paper was done
during a visit to Universit\'e Paris 7 which was supported  by Centre
franco-indien pour la promotion de la recherche avanc{\'e}e (CEFIPRA)
under Project 2501-1 {\it Algebraic Groups in Arithmetic and Geometry}.\newline
}}

\date{}

\maketitle

\begin{abstract}
We study deformation rings of $n$-dimensional mod
$\ell$ representations $\rho$ of the arithmetic fundamental group
$\pi_1(X)$ where $X$ is  a geometrically irreducible, smooth curve
over a finite field $k$ of characteristic $p$ ($ \neq \ell$).
We are able to show in many cases that the resulting rings are finite
flat over $\BZ_\ell$. The proof principally uses a lifting result of
the authors in Part~I of this two-part work, Taylor-Wiles systems
and the result of Lafforgue. This implies a conjecture of A.J. ~de
Jong in \cite{dejong} for representations with coefficients in power
series rings over finite fields of characteristic $\ell$,  that have
this mod $\ell$ representation as their reduction.   In fact a proof of all cases
of the conjecture for 
$\ell>2$ follows from a result announced 
(conditionally) by Gaitsgory in \cite{gaitsgory}.  The methods are completely different.
\end{abstract}

\tableofcontents

\section{Introduction}

Let $X$ be a geometrically irreducible, smooth curve over a finite
field $k$ of characteristic $p$ and cardinality $q$. Denote by 
$K$ its function field and by $\widetilde{X}$ its smooth
compactification and set $S:=\widetilde X\setminus X$. 
Let $\pi_1(X)$ denote the arithmetic fundamental group of $X$. 
Thus $\pi_1(X)$ sits in the exact sequence $$0 \rightarrow
\pi_1(\overline{X}) \rightarrow \pi_1(X) \rightarrow G_k \rightarrow 0,$$
where $\overline{X}$ is the base change of $X$ to an algebraic closure
of $k$, and $G_F$ denotes the absolute Galois group of any
field~$F$. In this paper we study deformation rings of  mod $\ell$
representations of $\pi_1(X)$, i.e., continuous (absolutely)
irreducible representations $\rhobar:\pi_1(X) \rightarrow
\GL_n({\BF})$ with ${\BF}$ a finite field of characteristic $\ell \neq
p$. We begin with the following conjecture of A.\ J.\ de Jong:

\begin{conj}[{\cite{dejong}}, Conj.~1.1]\label{DeJConj} 
Let $\rho\!:\pi_1(X)\to
\GL_n(\BF[[x]])$ be a continuous representation with residual
representation $\rhobar$. Then $\rho(\pi_1(\overline{X}))$ is finite.
\end{conj}

\begin{remark} \label{rOnFinRamif} 
{\em \begin{enumerate}\item In \cite{dejong} de Jong proves the above
    for $n\leq2$ by extending Drinfeld's reciprocity theorem in
    \cite{drinfeld} to $\BF((x))$-coefficients. 
\item It is an important feature, observed in Lemma 2.12 of
  \cite{dejong}, of the representations considered in Conjecture
  \ref{DeJConj} that the image of any inertia group (for a place in
  $S$) is finite. We will exploit this several times.\end{enumerate}}
\end{remark}

To state a reformulation of the above, we  we need some notation. Let $\CO$ be
the ring of integers of a finite extension of the fraction field of
$W(\BF)$ inside $\overline\BQ_\ell$, let $\rhobar$ be as above and fix a
lift $\eta\!:\pi_1(X)\to \CO^*$ of finite order of the $1$-dimensional
representation $\det\rhobar$. Then in \cite{dejong}, following
\cite{mazur}, it is  explained how to attach a deformation ring
$R^\eta_{X,\CO}(\rhobar)$, or simply $R^\eta_X$ for deformations of
$\rhobar$ of determinant $\eta$ and defined on~$\pi_1(X)$. 
In \cite{dejong} the following is shown:
\begin{theorem}[{\cite{dejong}}]\label{ConjEqConj}
Suppose $\rhobar$ is absolutely irreducible when restricted to
$\pi_1(\overline X)$. Then Conjecture~\ref{DeJConj} is equivalent 
to $R_X^\eta$ being finite (as a module) over~$\BZ_\ell$.
\end{theorem}
The theorem combined with the result quoted in
Remark~\ref{rOnFinRamif}~(a) shows: 
\begin{cor}[{\cite{dejong}}]\label{cConjForNeqTwo}
Suppose $n=2$ and $\rhobar$ is absolutely irreducible when restricted to
$\pi_1(\overline X)$. Then $R_X^\eta$ is finite over~$\BZ_\ell$.
\end{cor}

\begin{remark}\label{FinImpliesFlatCI}{\em 
Using obstruction theory, de Jong shows in \cite{dejong} that if
$R_X^\eta$ is finite over $\BZ_\ell$, then it is also flat over
$\BZ_\ell$ and a complete intersection.} 
\end{remark}

\subsection{Results}

In Theorem~\ref{OnDeJConj} we prove in many new instances that
the ring $R_X^\eta$ is finite over $\BZ_\ell$. To avoid some
technicalities, here we only state the following special case:
\begin{theorem}\label{SpecialDeJConj}
Let $\rhobar:\pi_1(X) \rightarrow \SL_n({\BF})$ be a
representation with ${\BF}$ a finite field of characteristic
$\ell \neq p$. Assume that
\begin{itemize}
\item[(i)] $\rhobar$ has full image, $\ell\notdiv n$, $|{\BF}|\ge4$,
  and $|{\BF}|>5$ if $n=2$, 
\item[(ii)] there exists $v\in S$ such that the
restriction of $\rhobar$ to a decomposition group at $v$ is
absolutely irreducible and the ramification is tame.
\end{itemize}
Then the ring $R^\eta_X$ is finite over $\BZ_{\ell}$, and in
particular Conjecture~\ref{DeJConj} holds for all $\rho$ with
reduction~$\rhobar$. 
\end{theorem}

\begin{cor}
Let $\rhobar$ be as in the previous theorem, and assume further that
at any place in $S$ ramification is either tame or of order prime to
$\ell$. Then $\rhobar$ has an $\ell$-adic lift
whose conductor is the same as that of $\rhobar$ and which is
automorphic, i.e., arises via reduction from a cuspidal eigenform of
the same conductor. 
\end{cor}

\begin{Proofof}{the corollary} 
Let $R_X^{0,\eta}$ denote the quotient of $R_X^\eta$ which
parameterizes deformations which are minimal at the places in
$S$. (This is a purely Galois-theoretic requirement. The conditions on
ramification are needed so that we can formulate, using
\cite{boecklekhare}, Props.~2.15 and~2.16, a minimality condition.) Using
Poitou-Tate and some obstruction theoretic arguments due to Mazur, it
is by now standard to show that $R_X^{0,\eta}$ has a presentation
$W(\BF)[[x_1,\ldots,x_n]]/(y_1,\ldots,y_n)$, where some of the $y_i$
could be zero. Because $R_X^\eta/(\ell)$ is finite, the same holds for
$R_X^{0,\eta}/(\ell)$. From this one deduces easily that $R_X^{0,\eta}$ must
be finite flat over $\BZ_\ell$. This proves the corollary when
combined with the results of \cite{lafforgue}.
\end{Proofof}

\begin{remark}{\em The corollary combined with finiteness theorems of
  Harder about dimensions of cusp forms with bounded conductor and
  fixed central character has consequences for conjectures in
  \cite{khare2}, \cite{moon}, \cite{moontaguchi}; see \cite{bk-moontag}}.
\end{remark} 

\subsection{On the proof, our hypotheses and related works}

\noindent{\bf On the proof of Theorem \ref{SpecialDeJConj}:}

The form of Conjecture~\ref{DeJConj}, combined with
Remark~\ref{rOnFinRamif}~(b) and Theorem
\ref{ConjEqConj}, lends itself to proving it not
necessarily over $X$, but over a suitable finite cover $Y$ of it,
i.e.\ to applying base change techniques. We repeatedly make use of this. 
In a first reduction, base change allows us to pass to a situation
where assumption (ii) of Theorem~\ref{SpecialDeJConj} holds at
all places of $S$  (and furthermore that~$|S|\ge 2$).
The results in Part~I of this work, cf. \cite{boecklekhare}, then
yield an $\ell$-adic lifting~$\rho$.

In a second reduction step, we apply a level lowering technique of
Skinner-Wiles in \cite{skinnerwiles}. It relies on an important 
principle that was discovered by Carayol in \cite{carayol} to switch
types of automorphic representations that give rise to a given
$\rhobar$. To apply Carayol's principle, we emply a Jacquet-Langlands
correspondence, so that we can work in an anisotropic setting. The
technique of Skinner-Wiles yields, after finite base change, a
minimal lift of $\rhobar$ over some finite cover $Y$ of~$X$. (As the
results in in \cite{lafforgue} provide us with all base change results
one expects, this does not require solvable base change.) 

Thus it suffices to prove the conjecture over this $Y$ and prove it
only for minimal deformation rings. This is again a significant
simplification as we need for instance no level raising results which
are in any case still not available, at least in any generality (as
far as we know!).

In a final step we construct Taylor-Wiles systems for
$\rhobar|_{\pi_1(Y)}$ using the Galois cohomology techniques of
Chapter~2 and automorphic methods of Chapter~3. (Again we use a
Jacquet-Langlands correspondence, so that we may consider
Hecke-modules in an anisotropic setting.) This allows one to prove
that a minimal deformation ring $R^{0}_Y$ for
$\rhobar|_{\pi_1(Y)}$ is finite over ${\bf Z}_{\ell}$. 
By what we have said this proves Theorem \ref{SpecialDeJConj}. 

Our techniques follow closely the original method of Wiles and Taylor
in \cite{wiles} and \cite{taylorwiles}, and later developments
\cite{diamond}, \cite{fujiwara}, \cite{skinnerwiles},
\cite{harristaylor}, which we have to generalize to our context.
In fact most of the work of this paper
is devoted to carrying out these generalizations. There is a small
modification needed to handle problems arising from ``torsion'' which
may be of relevance even in the number field case: in an appendix we
explain this innovation in the context considered in
\cite{wiles} and \cite{taylorwiles}. 

\bigskip

\noindent{\bf Remarks on the conditions on $\rhobar$ in
  Theorem~\ref{SpecialDeJConj}:} 
1. Condition~(i) is needed to apply our results from Part~I to obtain
  a lifting of $\rhobar$ to an $\ell$-adic Galois
  representation. Basically one could weaken it to conditions
  (a)--(c) of \cite{boecklekhare}, Thm.~2.1. Condition (d) of
  \cite{boecklekhare}, Thm.~2.1, is superfluous, since we may apply
  base change to achieve it.

\smallskip

2. Condition~(ii) is needed to be able to apply the Jacquet-Langlands
correspondence established in Badulescu's work,
cf. \cite{badulescux}, \cite{badulescuy}. The conditions imposed by
Badulescu require that there is at least one place at which the
representation is supercuspidal. So, in principle Badulescu's work
should allow one to prove Theorem~\ref{SpecialDeJConj} under the weaker hypothesis that
the restriction of $\rhobar$ to some decomposition group is absolutely
irreducible. (There are problems in this generality with our lifting
results from Part~I. But it seems likely that they could be resolved.)
Our main reason for working with condition~(ii) was, that it allows a
simple explicit description of the local types over $\GL_n$ and a
corresponding division algebra under the JL-correspondence. 

More generally, if the complete Jacquet-Langlands correspondence was known for
function fields, that would determine completely which automorphic
forms transfer from $GL_n$ to a given anisotropic inner form, we could
presumably weaken hypothesis (ii) considerably. 
For instance we would have then been able to cover cases where
the local ramification at a place was tame and such that a topological
generator of tame inertia was mapped to a regular unipotent element.
It is conceivable 
(the authors are far from being experts in this!) that the methods of 
Chapter VI.1 of \cite{ht} can also be adapted to the function field 
context to prove more instances of the JL correspondence.

\smallskip

3. Our methods do seem to require that one is able to define the
``minimal Hecke algebra'' using anisotropic groups which makes
cuspidal automorphic functions with given level and central character
behave like functions on finite sets. We have tried and failed to
avoid this by exploiting the fact that in the function field
situation, (cuspidal) automorphic forms with arbitrary coefficients
(in which $p$ is invertible) make sense even for split groups like
$GL_n$, and the result of Harder that cusp forms on $GL_n({\bf A}_K)$
(with $K$ the function field of $X$ and ${\bf A}_K$ its adeles) 
have compact support mod center and $GL_n(K)$. The main technical
difficulty is that we are unable to prove that a mod $\ell$ cusp form
lifts to an $\ell$-adic cusp form at least up to ``Eisenstein error'':
this is explained more precisely in the sentence after Proposition
\ref{PropRedOfCuspForms}.

\bigskip

\noindent{\bf Comparison to other works:}

D. Gaitsgory has informed us that he can prove Conjecture~\ref{DeJConj} for
arbitrary $\rhobar$ (i.e., not even necessarily irreducible) 
provided $\ell>2n$ by methods which are in a sense akin to de Jong's proof
in the case $n=2$, i.e., they directly prove automorphy of the representations 
considered in \ref{DeJConj}. As far as the authors know, this has not
as yet been written up.

Our results do not have  a restriction like 
$\ell>2n$. On the other hand we do have somewhat unpleasant and serious restrictions on
$\rhobar$. The methods of Gaitsgory and this paper are presumably completely different. In a sense we think
the methods in de Jong's proof of the conjecture for $n=2$ and, and also Gaitsgory's work, are more natural 
(given the geometric context of function fields that we are in) than ours! 
Nevertheless our results illustrate the point that the methods in
\cite{wiles}, \cite{taylorwiles} as refined in \cite{diamond}, \cite{fujiwara}, \cite{skinnerwiles},
are flexible enough to be adapted to a variety of situations, for instance the present context of function fields.

\smallskip 

{\it This was the situation at the time when a version of this work was submitted. Since then a proof of all
cases of de Jong's conjecture has been announced in \cite{gaitsgory} for $\ell>2$, which depends on some
results not yet written up. The methods are completely different to ours, and our proof works in many cases 
even when $\ell=2$.}

\smallskip

We note that there are 2 separate  works of Genestier and Tilouine 
(for $GSp_4$), and Clozel, Harris and Taylor (for $GL_n$), on 
generalizations of \cite{wiles} and \cite{taylorwiles} to $n$-dimensional
representations of absolute Galois groups of totally real number
fields. In fact after a rough version of this paper was written
\cite{harristaylor} became available to us (thanks to Michael Harris),
and its a comfort to us that we may quote verbatim from Section II of
{\it loc.\ cit.}\ for the Hecke action at places that are introduced
in building Taylor-Wiles systems. One of the key technical differences
between the present work and \cite{harristaylor} when building
Taylor-Wiles systems is that we allow for the possibility that the
$\ell$th roots of 1 are in the base field which entails slight
adjustments on the Galois and automorphic side. 

The key {\it qualitative} difference between the mentioned works and
ours is that we can prove automorphy of residual representations like
$\rhobar$ in the theorem, while in the other works this has to be at
the moment an imprtant assumption that seems extremely difficult to
verify in their number field case; further we are mainly interested in
establishing algebraic properties of deformation rings, while in the
number field case these are established {\it en route} to 
proving modularity of $\ell$-adic representations (which is known in
our context by \cite{lafforgue}!). Thus our uses of the methods
pioneered by Wiles can be deemed to a certain extent to be warped! 

\medskip

{\em Throughout this paper we use the notation of Part I of this work without further mention.}

\section{Galois cohomology}

\subsection{On deformations of the determinant}

Unlike in part one of this work, we will no longer work with fixed
determinants when considering deformation problems. At the same time,
we want to keep finite the number of such determinants.
The standard way over function fields to achieve this is to require
that under any deformation of the determinant a certain a priori
chosen place is totally split. The choice of place is the content of
the following lemma: 
\begin{lemma}\label{lOnPlaceW}
There exist infinitely many places $w\in X$ whose residue
field $k_w$ satisfies $\ell\notdiv[k_w:k]$ and $k(\zeta_\ell)\cap k_w=k$.
\end{lemma}
\begin{Proof}
Let $K'\subset K(\zeta_{\ell^2})$ denote the unique constant field
extension of $K$ of degree $\ell$. Then $K'$ and $K(\zeta_\ell)$ are 
Galois over $K$. Let $\sigma$ be any element of
$\Gal(K(\zeta_{\ell^2})/K)$ which lies in neither of the normal
subgroups $\Gal(K(\zeta_{\ell^2})/K')$ or
$\Gal(K(\zeta_{\ell^2})/K(\zeta_\ell))$. By the \v{C}ebotarov density
theorem there exist infinitely many places $w\in X$ whose Frobenius
automorphism maps to $\sigma$. Any such will have the desired
properties.
\end{Proof}

For a place $w$ as above one and $\overline I_w:=I_w/([G_w,I_w]I_w^\ell)$ one has a
split exact sequence
\begin{equation}\label{SESforGabW}
0\longto \overline I_w \longto G_w^\ab\longto \hat\BZ\longto 0
\end{equation}
where the group $\overline I_w\cong k_w^*/k_w^{*\ell}$ is a quotient of
$\BZ/(\ell)$ and trivial unless $\zeta_\ell\in K$. We fix a splitting
$s_w\!:G_w^\ab\to \overline I_w$.
\begin{lemma}\label{CorOnLandLdualForDet}
Let $w, s_w$ be as above. Let $\rhobar\!:G_w \to \{1\}\in\BF$ denote
the trivial character and $\rho^w\!:G_w\to \GL_1(R^w)$ the universal
deformation of $\rhobar$ for deformations which factors via the
splitting $s_w$. Denote by $L_{w}\subset H^1(G_w,\BF)$ the subspace
corresponding to the tangent space of the universal deformation. 

Then $\dim L_w=\dim_\BF\BF(1)^{\pi_1(X)}$, $L_w$ is disjoint from
$H^1_\unr(G_w,\BF)$, and $L_w^\perp\subset H^1(G_w,\BF(1))$ is of dimension one and disjoint
from $H^1_\unr(G_w,\BF(1))$.
\end{lemma}
\begin{Proof}
From the definition of $R^w$ it follows that $L_w$ is a complementary
sub vector space for $H_\unr^1(G_w,\BF)$ in $H^1(G_w,\BF)$. Since
local Tate-duality is perfect, $L_w^\perp$ must be a complementary
vector space for $H_\unr^1(G_w,\BF(1))$ in $H^1(G_w,\BF(1))$. 

By a standard spectral sequence argument as used in
\cite{boecklekhare}, sequence (4), p.~13, one obtains a short exact
sequence
$$0\longto \BF_{G_w}\cong H^1_{\unr}(G_w,\BF)\longto H^1(G_w,\BF)
\longto \BF(-1)^{G_w}\longto 0$$
and similarly for $\BF(1)$ instead of $\BF$. From this one deduces
$L_w\cong \BF(-1)^{G_w}$ and $L_w^\perp\cong \BF^{G_w}\cong\BF$. All
assertions are now obvious.
\end{Proof}

\smallskip
 
{\em From now on, for the remainder of this article, we fix  a place
$w\in X$ as in Lemma~\ref{lOnPlaceW} and a splitting $s_w$ is in
(\ref{SESforGabW}).}

\subsection{On Taylor-Wiles' auxiliary primes}\label{TWPrimes}

In Section~\ref{TWSyst}, we will construct Taylor-Wiles systems,
\cite{taylorwiles}, in the minimal case. As such they consist of a
Galois-theoretic and a Hecke part. The current section provides the
Galois-theoretic tools needed. Recall that in \cite{boecklekhare} the
extension $E/K$ is defined as the splitting field~$\rhobar$.
\smallskip

We begin with the following Lemma:
\begin{lemma}\label{LemmaOnTWPrimes}
Let $v$ be a place such that $q_v\equiv 1\!\!\pmod\ell$, $\rhobar$ is
unramified at $v$ and $\rhobar(\Frob_v)$ has distinct eigenvalues
which are all contained in $\BF$. Let $(R,\Fm)$ be in $\CA$ and
let $\rho_v\!:G_v\to \GL_n(R)$ be a lift of $\rhobar_v$. Then up to
strict equivalence the image of $\rho_v$ is diagonalizable. 
\end{lemma}
Recall that $\rho_v,\rho_v'\!:G_v\to\GL_n(R)$ are strictly equivalent
if there exists $M\in\GL_n(R)$ congruent to the identity modulo
$\Fm$ such that $M\rho_vM^{-1}=\rho'_v$. 

Note also that $\rho_v$ will factor through the tame quotient $\overline
G_{q_v}$ of $G_v$ since $\rhobar$ is unramified and the kernel of
$\pi\!:\GL_n(R)\to\GL_n(\BF)$ is prime to $p$. The lemma implies that
in fact $\rho_v$ factors through the abelianization
$\hat\BZ\times\BZ/(q_v-1)$ of $\overline G_{q_v}$. 

\begin{Proof}
Let us assume that $\rhobar_v$ takes its image in the diagonal
matrices. This shows in particular, that the exponent $e$ of the cyclic
group $\rhobar(G_v)$ is prime to $\ell$. Because the kernel of $\pi$,
defined above, is a pro-$\ell$ group, the representation $\rho_v\!:\overline
G_{q_v}\to \GL_n(R)$ must factor via the quotient
$G:=\BZ_\ell(1)\rtimes (\BZ_\ell\times\BZ/(e))$ of $\overline G_{q_v}$. 

Let $\sigma$ be a generator of $\BZ/(e)$ and
$s$ of $\BZ_\ell\times\BZ/(e)$. We may regard $\sigma$ as well as
$s$ as elements of $G$. Because $\rhobar(\sigma)$ has distinct
eigenvalues and $\rho_v(\sigma)$ has finite order $e$, using
strict equivalence we may assume that $\rho_v(\sigma)$ is diagonal.
Since $\rhobar(\sigma)$ has distinct eigenvalues, the same holds for
$\rho_v(\sigma)$. But this implies that $\rho_v(s)$ is diagonal as
well, because $\rho_v(s)$ commutes with~$\rho_v(\sigma)$.

We now claim by induction on $i\geq1$ that $\rho_v(t)\pmod{\FM^i}$ is
diagonal, the assertion being trivially true for $i=1$. So let us
assume that the assertion holds for $i$ and write $\rho_v(t)=D+B$,
where $D$ is a diagonal matrix and $B$ is zero along the diagonal and
has entries in $\FM^{i}$. Because $B^2\equiv0\!\!\pmod{\FM^{i+1}}$,
the relation $sts^{-1}=t^{q_v}$ yields 
$$\rho_v(s)(D+B)\rho_v(s)^{-1}\equiv
D^{q_v}+\sum_{i=0}^{q_v-1}D^iBD^{q_v-i-1}\!\!\!\pmod{\FM^{i+1}}.$$ 
As $D\equiv I\!\!\pmod\FM$, and $q_v-1\in\FM$ the right hand side
is congruent to $D^{q_v}+B$ modulo $\FM^{i+1}$. Comparing off-diagonal
entries, we see that $B\pmod  {\FM^{i+1}}$ commutes with
$\rho_v(s)\pmod  {\FM^{i+1}}$. This shows that $B\pmod  {\FM^{i+1}}$
is diagonal, and hence zero.
\end{Proof}

For a place $v$ as in the previous lemma, 
an integer $m\in\BN$ such that $\ell^m|(q_v-1)$, and
an eigenvalue $\lambda$ of $\rhobar(\Frob_v)$, we define the collection of
lifts $\CC_{v,\lambda,m}$ of $\rhobar_{v}$ as follows: 

Let $s,t\in G_v^\ab\cong\BZ/(q_v-1)\times\hat\BZ$
be such that $t$ generates inertia and $s$ maps to $\Frob_v$ in $G_v/I_v$.
Denote by $\lambda_1,\ldots,\lambda_n$ Teichm\"uller lifts of the
eigenvalues if $\rhobar(\Frob_v)$ such that
$\lambda\equiv\lambda_1\pmod\ell$. Set
$$R_{v,m}:=W(\BF)[[x_1,\ldots,x_n,y]]/\Big(
(1+y)^{\ell^m}-1\Big),$$
and define $\rho_{v,\lambda,m}\!:G_v\onto G_v^\ab
\longto \GL_n(R_{v,m})$ by
$$s\!\mapsto \Diag(\lambda_1(1+x_1),\ldots,\lambda_n(1+x_n)),\,
t\!\mapsto \Diag((1+y),1,1,\ldots,1).$$
Finally define $\CC_{v,\lambda,m}\!:\CA\to\Sets$ by   
\begin{eqnarray*} R&\mapsto& \CC_{v,\lambda,m}(R):=\left\{\rho\!:G_v\to\GL_n(R)\,|\,
\exists\alpha\in\Hom_\CA(R_{v,m},R),\right.\\
&&\left.\phantom{\CC_{v,\lambda,m}(R):=\left\{\rho\!:\right.} \exists M\in1+M_n(\FM_R):\rho=M(\alpha\circ\rho_{v,\lambda,m})M^{-1}\right\}.
\end{eqnarray*}
We define $L_{v,\lambda,m}\subset H^1(G_v,\Ad)$ as the
subspace spanned by the $1$-cocycles 
\begin{equation}\label{DefOfLvLambda}
\left\{c:G_v\to\Ad: g\mapsto 1/\eps(\rho(g)\rho_0^{-1}(g)-I)
| \rho\in\CC_{v,\lambda,m}(\BF[\eps]/(\eps^2)) \right\},\end{equation}
where $\rho_0$ is the tautological lift induced from the splitting
$\BF\to\BF[\eps]/(\eps^2)$. It is easy to see that $\dim
L_{v,\lambda,m}=\dim H_\unr^1(G_v,\Ad)+1=n+1$. Moreover
$L_{v,\lambda,m}$ is independent of $m$ as long as $m>0$. Finally,
note also that \cite{boecklekhare},~Lem.~2.13, yields: 
\begin{lemma}\label{DependsNotOnV}
Let $\sigma\in\Gal(E(\zeta_\ell)/K)$ be the image of $\Frob_v$. Then the subspace 
$$L_{v,\lambda,m}^\perp\subset H^1_\unr(G_v,\Ad(1))\cong
\Ad/(\Frob_v-1)\Ad$$
of codimension one only depends on $\sigma$ and the choice of
$\lambda$ (among the eigenvalues of $\rhobar(\sigma)$).
\end{lemma}
Because of the above lemma, we also write $L^\perp_{\sigma,\lambda}$
or $L^\perp_{v,\lambda}$ for $L_{v,\lambda,m}^\perp$. Note that since the
cyclotomic character $\chi$ is trivial on $G_v$, the restrictions of
$\Ad$ and $\Ad(1)$ to $G_v$ agree.

\medskip

The central result of this section is the following which is modeled
at \cite{harristaylor}, Thm.~IV.5.3:
\begin{lemma}\label{LemONTWprimes} 
Let $\rhobar:\pi_1(X) \rightarrow \GL_n({\BF})$ be a continuous
representation and $(\CC_v,L_v)_{v\in S\cup \{w\}}$ a set of
deformation conditions. Suppose that
\begin{enumerate}
\item \label{AssOnea} $L_v=H^1_\unr(G_v,\Ad)$ for all $v\in S$ (a minimality condition
  on the local deformations).
\item \label{AssOneb} $L_w=H^1_\unr(G_v,\ad)\oplus L_{w,d}$, where
  $L_{w,d}$ is the subspace of $H^1(G_w,\BF)$ considered in
  Lemma~\ref{CorOnLandLdualForDet}, and $\BF\subset\Ad$ via the
  diagonal embedding.
\item \label{AssTwo} For any $\pi_1(X)$-subrepresentation $V$ of
  $\Ad$, there exists a regular semisimple
  $g_V\in\rhobar(\pi_1(\overline X))$ such that $V^{g_V}\neq0$.
\item\label{AssThreeA} If $\zeta_\ell\in E$, then
  $
  H^1(\Gal(E/K(\zeta_\ell)),\ad)=0$.
\item If $\zeta_\ell\in K$, then\label{AssThreeB} $\ad$ has no
 $1$-dimensional subrepresentation.
\item \label{AssFour} The image of $\rhobar$ has no quotient of order $\ell$.
\end{enumerate}
Then for any given $m\in\BN$ there exists a set $Q_m$ of $\dim
H^1_{\{L_v\}}(\{w\},\Ad)$ places of $X\setminus\{w\}$ such that 
\begin{itemize}
\item[(a)] $q_{v}\equiv 1\!\!\pmod{\ell^m}$ for all $v\in Q_m$,
\item[(b)] $\rhobar(\Frob_{v})$ has distinct eigenvalues for each $v\in Q_m$, and 
\item[(c)] $H^1_{\{L_{v}^\perp\}}(Q_m,\Ad(1))=0$ where
  $L_v=L_{v,\lambda_{v}}$ for each $v\in Q$ and $\lambda_v$ is some
  eigenvalue of $\rhobar(\Frob_{v})$. 
\end{itemize}
Moreover $\dim H^1_{\{L_v\}}(\{w\},\Ad)=\dim
H^1_{\{L_{v}\}}(Q_m\cup\{w\},\Ad)$.
\end{lemma}
The proof in fact shows that the sets $Q_m$ above may be chosen disjoint
from any given finite set of places $S'\supset S$.

\begin{Proof}
The proof follows closely the analogous proof given in
\cite{harristaylor} which in turn is similar to that in
\cite{taylorwiles}.
Note first that  $L_v=H_\unr(G_v,\Ad)$ implies that $\dim
L_v=h^0(G_v,\Ad)$. Because of Lemma~\ref{CorOnLandLdualForDet}, by
\cite{boecklekhare},~Rem.~2.6, we have
$$\dim H^1_{\{L_v\}}(\{w\},\Ad)=\dim
H^1_{\{L^\perp_v\}}(\{w\},\Ad(1)).$$ 

Define $E_m:=E(\zeta_{\ell^m})$ and let $Y_m\to X$ be the
corresponding Galois cover. We first claim that the composite
$$H^1_{\{L^\perp_v\}}\!(\emptyset,\!\Ad(1)) \!\into \!H^1\!(\pi_1(X),\!\Ad(1))\!\to\!
H^1\!(\pi_1(Y_m),\!\Ad(1))^{\Gal(E_m/K)}$$ 
is injective, where the second morphism is restriction.

For any $v\in S\cup \{w\}$ we have $L_v=L_{v,0}\oplus L_{v,d}$
with $L_{v,0}:=L_v\cap H^1(G_v,\ad)$ and $L_{v,d}:=L_v\cap
H^1(G_v,\BF)$, and so we may prove the claim separately for the
subrepresentations $\BF$ and $\ad$ of $\Ad$. We first consider~$\ad$.

Condition~\ref{AssThreeA} yields $H^1(\Gal(E_1/K),\ad(1)))=0$, as can
be seen by applying for instance by \cite{boeckle},
Prop.~1.8~(i),(ii). This proves the claim for $m=1$ and $\ad$.
%
For $m>1$, inflation-restriction and taking invariants yields the
left-exact sequence
$$0\to H^1(\Gal(E_m/E_1),\ad(1))^{\Gal(E_1/K)}\to
H^1(\pi_1(Y_1),\ad(1))^{\Gal(E_1/K)}$$
$$\to H^1(\pi_1(Y_m),\ad(1))^{\Gal(E_m/K)}. $$
We will show that the left hand term vanishes. Since
$K(\zeta_{\ell^m})$ is Galois over $K$, the 
group $\Gal(E_m/E_1)$ lies in the center of~$\Gal(E_m/K)$.
Moreover by the definition of $E_1$, the
action of $\Gal(E_m/E_1)$ on $\Ad(1)$ is trivial. Therefore 
we find
\begin{eqnarray*}\lefteqn{H^1(\Gal(E_m/E_1),\ad(1))^{\Gal(E_1/K)}}\\
&=&H^1(\Gal(E_m/E_1),\BF)\otimes_\BF(\ad(1))^{\Gal(E_1/K)}.
\end{eqnarray*}
By \ref{AssThreeB} the last expression is zero. This proves the claim
for the $\ad$-component.

\smallskip

We will now consider the diagonal $\BF$-component. By
inflation-restriction we need to show
$$0=H^1_{\{L^\perp_{v,d}\}}(\{w\},\BF(1))\cap H^1(\Gal(E_m/K),\BF(1))$$
inside $H^1(\pi_1(X\setminus\{w\}),\BF(1))$. Again by
inflation-restriction, the second term allows the isomorphism
$$H^1(\Gal(E_m/K),\BF(1))\cong
(H^1(\Gal(E_m/K(\zeta_\ell)),\BF)\otimes\BF(1))^{\Gal(K(\zeta_\ell/K))}.$$ 
By assumption \ref{AssFour} the right  hand side is isomorphic to
$$H^1(\Gal(K(\zeta_{\ell^m})/K(\zeta_\ell)),\BF)
\otimes(\BF(1)^{\Gal(K(\zeta_\ell)/K)}).$$
If $\zeta_\ell\notin K$, the proof of the claim is thus complete. In
the case $\zeta_\ell\in K$, a non-zero class in $ H^1(\Gal(E_m/K),\BF(1))$
describes a non-zero character
$\Gal(K(\zeta_{\ell^m})/K(\zeta_\ell))\onto\BF$. Because
$\ell\notdiv[k_w:k]$, such a class restricts to a non-zero class in
$H^1_\unr(G_w,\BF(1))$. By Lemma~\ref{CorOnLandLdualForDet} it maps to a
nonzero class in $H^1(G_w,\BF(1))/L^\perp_{w,d}$ and so such a class
cannot lie in $H^1_{\{L_{v,d}^\perp\}}(\{w\},\BF(1))$. This completes
the proof of the claim.

\medskip

Let now $\psi$ be a $1$-cocycle whose class $[\psi]\in
H^1_{\{L_v^\perp\}}(\emptyset,\Ad(1))$ is non-zero. By the above
claim, the restriction of $\psi$ to $\pi_1(Y_m)$ is non-trivial. Since
$\pi_1(Y_m)$ acts trivially on $\Ad(1)$, the class $[\psi]$ induces a
non-trivial $\Gal(E_m/K)$-equivariant homomorphism $\pi_1(Y_m)\to\Ad(1)$. 
Let $E_\psi$ denote the fixed field of its kernel, and $V(1)$ its
image in $\Ad(1)$. Then the induced morphism
$\Gal(E_\psi/E_m)\to V(1)$ is bijective and non-zero.

Choose a regular semisimple $g\in \rhobar(\overline X)$ such that
$V^g\neq0$. Because $\rhobar$ has finite image, we have
$g=\rhobar(\sigma)$ for some $\sigma\in\Gal(E_1/K(\zeta_\ell))$.
For such a $g$ we now consider
$$V(1)_g\cong V_g \longinto \Ad(1)_g\cong
\Ad_g=H_\unr^1(G_v,\Ad(a)),$$
where $v$ is any place with $\Frob_v\mapsto \sigma$
($\Frob_v\in\pi_1(X\setminus\{w\})$). For an eigenvalue 
$\lambda$ of $g$ we denote by $L_{g,\lambda}^\perp\subset
\Ad(1)_g$ the corresponding subspace as defined in
(\ref{DefOfLvLambda}). One easily shows $\cap_\lambda
L_{g,\lambda}^\perp=0$, where the intersection ranges over all
eigenvalues of $g$, e.g., by proving the dual  assertion.
We claim that there exists $x_V\in V_g$ and an
eigenvalue $\lambda$ of $g$ such that
$$\psi(\Frob_v) + x_V \notin L^\perp_{g,\lambda} \subset H^1_\unr(G_v,\Ad(1)).$$

Assume otherwise. Then $\psi(\Frob_v)+V_g\subset L_{g,\lambda}^\perp$
for all eigenvalues $\lambda$. This implies $V_g\subset L_{g,\lambda}^\perp$
for all $\lambda$, and therefore $V_g\subset \cap_\lambda
L_{g,\lambda}^\perp=0$, contradicting $V_g\cong V^g\neq0$.

\smallskip

Let $\tau\in\pi_1(X\setminus\{w\})$ by any element which acts
trivially on $E_m$ and maps to $x_V$ under $\psi$. Because $E_\psi/K$
is Galois, by the \v{C}ebotarov theorem there exists a place $v'\in
X\setminus\{w\}$ such that the image of $\Frob_{v'}$ in
$\Gal(E_\psi/K)$ agrees with that of $\tau\Frob_v$. By the above
$[\psi]\notin L_{v',\lambda}^\perp$. A simple inductive argument now
finishes the proof of the lemma.
\end{Proof}

To complement the above lemma, in the particular case in which we are
interested, we also prove the following:
\begin{lemma}\label{SL_nAdmitsTWclasses}
Suppose $\rhobar\!:\pi_1(X)\to \SL_n(\BF)$ is surjective, $\ell\notdiv
n$, $|\BF|\ge4$, and $|\BF|>5$ for $n=2$. Then
\begin{enumerate}
\item $\rhobar(\pi_1(\overline X))=\SL_n(\BF)$ contains a regular semisimple element.
\item $\Ad=\ad\oplus\BF$, $\ad$ is irreducible, and  for any regular
  semisimple $g\in\SL_n(\BF)$ one has $\BF^g=\BF$ and
  $M^0_n(\BF)^g\neq0$.
\item $H^1(\SL_n(\BF),M^0_n(\BF))=0$. 
\item $\SL_n(\BF)$ has no quotient isomorphic to $\BZ/(\ell)$.
\end{enumerate}
\end{lemma}
\begin{Proof}
Under our conditions on $n$ and $\BF$, the group $\SL_n(\BF)$ has no
abelian quotients, which shows (iv) and $\rhobar(\overline{X})=
\SL_n(\BF)$. Parts (i) and (ii) are clear and (iii) follows from \cite{cps}.
\end{Proof}

\subsection{Removing local ramification}\label{RemoveRam}

The current section lays the Galois-theoretic ground work for the base
change techniques which we will apply repeatedly in the proof of our
main result in the subsequent chapter. In this respect the following
proposition and corollary will be of much use  to us.
\begin{prop}\label{EnlargeBaseFieldProp} 
Let $R$ be in $\CA$, let $\rho\!:\pi_1(X)\to \GL_n(R)$ a continuous
representation and let $S_1\subset S$ be the set of places at which 
$\rho(I_v)$ is finite. Let $T_s,T_i\subset \widetilde X$ be finite and
disjoint and let $m$ be some positive integer. Then there exists a
finite (possibly ramified) Galois cover $Y\to X$, say with function
field $L$, such that
\begin{enumerate}
\item $L/K$ is totally split at all places above of $T_s$.
\item At places in $T_i$, the residue degree of the extension
  field $L/K$ is a multiple of~$m$.
\item The restriction $\rho_{|\pi_1(Y)}$ is unramified at the places
  above $S_1\setminus T_s$. 
\end{enumerate}
\end{prop}
Note that the assumptions are satisfied for any representation $\rho$
which is unramified outside a finite set of places, provided $R$
is finite or the ring of integers of a local field of positive
characteristic, cf.\ Remark~\ref{rOnFinRamif}~(b). Before giving the proof, we
state the following important corollary, which will be used in the sequel:
\begin{cor}\label{EnlargeBaseFieldCor} 
Suppose $R,\rho,S_1$ are as in the previous proposition. Let $T_r$ be
a set of places at which $\rho$ is ramified, $T_i\subset
\widetilde X$ any finite subset disjoint from $T_r$ and $m$ some
positive integer. Then there exists a finite (possibly ramified)
Galois covering $Y\to X$ with corresponding extension $L/K$ of
function fields such that
\begin{enumerate}
\item $Y$ is geometrically connected over $k$,
\item $\rho(\pi_1(X))=\rho(\pi_1(Y))$, 
  $\rhobar(\pi_1(\overline X))=\rhobar(\pi_1(\overline Y))$,
\item $\rho_{|\pi_1(Y)}$ is unramified above the places in
  $S_1\setminus T_r$.
\item For all places $v'$ in $L$ above a place $v\in T_r$ one has
  $\rho_{|G_{v'}}\cong \rho_{|G_v}$.
\item At places above $T_i$, the residue degree of $L/K$ is a multiple
  of~$m$.
\end{enumerate}
\end{cor}
\begin{Proofof}{Corollary~\ref{EnlargeBaseFieldCor}}
We claim that there is a finite set $T'\subset X$ disjoint from $T_i$
such that the elements $\rho(\Frob_v)$, $v\in T'$, topologically generate
$\rho(\pi_1(X))$. Since $\rhobar(\pi_1(X))$ is finite, by the
\v{C}ebotarev density theorem there exists a finite set $T''\subset X$
disjoint from $T_i$ with the
above property for $\rhobar$. Let $X'\to X$ be the smallest finite
Galois covering over which $\rhobar$ becomes trivial. Then
$\rho(\pi_1(X'))$ is a pro-$\ell$ group, and since the pro-$\ell$
completion of $\pi_1(X')$ is topologically finitely generated, so is
$\rho(\pi_1(X'))$. Therefore its Frattini  quotient is finite, and
again by the \v{C}ebotarev density theorem we may choose a finite set
$T'''\subset X$ disjoint from $T_i$ such that the elements
$\rho(\Frob_v)$, $v\in T'''$, lie in $\rho(\pi_1(X'))$ and span the Frattini
quotient. Therefore by Burnside's basis theorem these elements
topologically generate $\rho(\pi_1(X'))$. The claim follows with
$T':=T''\cup T'''$. 

Recall that $E$ is the splitting field of $\rhobar$.
Let now $T_s\supset T'\cup T_r$ be a finite set of places disjoint
from $T_i$ such that \begin{itemize}
\item[(a)] the greatest common divisor of the $q_v$, $v\in
T_s$, is $q$, and
\item[(b)] $T_s$ contains places $w_1,\ldots,w_j$ such that the
  $\Frob_{w_i}$, $i=1,\ldots,j$, lie in $G_E$ and such that the
  greeatest common divisor of the $q_{w_i}$, $i=1,\ldots,j$, agrees
  with $\# (E\cap\overline\BF_p)$. 
\end{itemize} 
Applying the above proposition with these data, the corollary
follows. (Note that condition (b) guarantees the second part of (ii).)
\end{Proofof}

We first prove the following lemma:
\begin{lemma}
Let ${\widetilde K}/K$ be finite with constant field $k$, let $u$ be a
place of ${\widetilde K}$ and ${\widetilde K}_u$ the corresponding
completion. Suppose $F/{\widetilde K}_u$ is Galois of prime degree
$e$. Then for any set of place $T$ of ${\widetilde K}$ not containing
$u$, there exists a Galois extension $L$ of ${\widetilde K}$ of degree
$e$ such that 
\begin{enumerate}
\item all places in $T$ are split,
\item there is a unique place $u'$ in $L$ above ${\widetilde K}$, and
\item the extension $L_{u'}/{\widetilde K}_u$ is isomorphic to $F/{\widetilde K}_u$. 
\end{enumerate}
\end{lemma}
\begin{Proof}
Let us first assume that $e=p$. Then $F/{\widetilde K}_u$ is an Artin-Schreier
extension, say given by an equation $x^p-x=a$ for some $a\in {\widetilde K}_u$.
Using the approximation theorem, we can find $b\in {\widetilde K}$ such that
$b-a$ lies in the maximal ideal of the local ring at $u$ and $b$ lies
in  the maximal ideal of the local rings of places in $T$. Then the
Artin-Schreier extension $L:={\widetilde K}(z)$ with $z^p-z=b$ has all
the required properties.

We now suppose $e\neq p$. In this situation, we quote the following
result from \cite{neukirch}, Ex.~6 (the required modifications for the
function field case are straightforward): For any finite set $T'$ of
places of ${\widetilde K}$, the morphism  
$$H^1(G_K,\BZ/(e))\longto \amalg_{v\in T'} H^1(G_v,\BZ/(e))$$ 
is surjective. Note that $H^1(G_K,\BZ/(e))=\Hom(G_K,\BZ/(e))$. 
Define $T':=T\cup\{u\}$ and set $\chi_v=0\in \Hom(G_v,\BZ/(e))$ for
$v\in T$ and choose for $\chi_u$ a character of $G_u$ whose fixed
field is $F$. Let $\chi$ be a global character which maps to the
$\chi_v$, $v\in T'$. Then it is easy to see that the fixed field of
$\chi$ satisfies all the conditions required from $L$.
\end{Proof}
\begin{Proofof}{Proposition~\ref{EnlargeBaseFieldProp}}
By possibly shrinking $S$, we may assume that $\rhobar$ is ramified at
all places of $S$. Let $v_1,\ldots,v_r$ denote the places in $S$ and
$v_1,\ldots,v_{r'}$, $r'\leq r$, those in $S_1\setminus T_s$. 
The groups $\rho(G_v)$ are all pro-solvable. By
repeated application of the above lemma, we first construct a finite
extension $L_1$ over $K$ which is totally split at all places in
$T_s$ and such that the local extension of $L_1/K$ above $v_1$, is
Galois with group isomorphic to $\rho(G_{v_1})$ under $\rho$. Then
construct $L_2/L_1$ which is again totally split at all places of $T_s$
and such that for a place $v'_2$ of $L_1$ above $v_2$ in $K$ the local
extension of $L_2/L_1$ above $v'_1$ is Galois with group isomorphic
to $\rho(G_{v_2})$ under $\rho$. One reaches inductively an
extension $L_{r'}$ which contains places $w_1,\ldots,w_{r'}$ above
$v_1,\ldots,v_{r'}$ such that the restriction of $\rho$ to each
$G_{w_i}$ is trivial. 

By another repeated application of the lemma, we may construct an
extension $L'_{r'}$ of $L_{r'}$, totally split at all places above $T_s$ and
such that the residue degree at places in $T_i$ grows by a
multiple of $m$. (Locally at places in $T_i$ one constructs the
unramified extensions of degree $m$.)
Then the Galois closure $L$ of $L'_{r'}$ above $K$ has all the
desired properties.
\end{Proofof}

\section{Automorphic methods}\label{AutoC}

We begin this chapter by briefly reviewing parts of the theory of
automorphic forms over function fields. In Sections~\ref{LocReps}
and~\ref{AutoBasics}, we  are concerned with the local, respectively
global theory. We present some cases of the local and global Langlands
correspondence and of a Jacquet-Langlands correspondence as proved by
Badulescu. First consequences of the above are presented in
Section~\ref{Hecke} on Hecke-algebras. We formulate in Section \ref{carprin}
 a principle of Carayol, \cite{carayol} in a form suitable for us (cf., the
2 lemmas in the section).
Section~\ref{TWSyst} explains how to carry over the method of
Taylor-Wiles systems and its later simplification, cf.\ 
\cite{taylorwiles}, \cite{diamond}, from the number to the function
field case, using crucially results of Section \ref{carprin}, assuming that we are in a situation where
we have an automorphic lift of $\rhobar$ of minimal level. Thereby we prove isomorphisms between certain universal
deformation rings and corresponding Hecke algebras, under this crucial assumption. We present a
slight technical improvement over the usual method that might also be
useful in the number field case. Namely we do not need to choose an
auxiliary level structure and therefore can also treat cases in which
the function field $K$ contains a primitive $\ell$-th root of unity (see also the appendix where this method is explained for the classical situation).
In Section~\ref{LevelLow}, we prove a theorem on `lowering the level'
for certain cuspidal Hecke eigenforms, which allows us to verify this assumption in a significant number of cases, enough to cover
the applications to Theorems \ref{OnDeJConj} and \ref{SpecialDeJConj} 
that are proved in the last section by pulling all the results of
this section together. 
Our method here is inspired by a 
base change idea of Skinner-Wiles, \cite{skinnerwiles}, and again relies
on results of Section~\ref{carprin}.

\medskip

We fix some notation to state the main theorem of the final section:
Let $\CO$ be a discrete valuation ring finite over $W(\BF)$ and with
maximal ideal $\FM$. For a finite subset $T$ of $X$, consider a
representation $\rho\!:\pi_1(X\setminus T)\to\GL_n(\CO)$ with residual
representation $\rhobar:=\rho\pmod{\FM}$. The extension $E/K$ is as in
Section~\ref{TWPrimes}. We say that
$\rho$ is {\em of type-$1$ at a place $v\in\widetilde X$}, if
{\em $\rhobar$ is unramified at $v$ and $\rho(I_v)$ is unipotent of rank $1$}.

In the final section, we prove the following
theorem and give some further applications of it.
\begin{theorem}\label{OnDeJConj}
Let $\CO,T,\rho,\rhobar$ be as above and suppose they satisfy
\begin{enumerate}
\item $\rho_v$ is of type-$1$ at all $v\in T$,
\item For any $\pi_1(X)$-subrepresentation $V$ of
  $\Ad$, there exists a regular semisimple
  $g_V\in\rhobar(\pi_1(\overline X))$ such that $V^{g_V}\neq0$.
\item If $\zeta_\ell\in K$, then $\ad$ has no $1$-dimensional
  subrepresentation,
\item If $\zeta_\ell\in E$, then $H^1(\Gal(E/K(\zeta_\ell),\ad)=0$.
\item The image of $\rhobar$ has no quotient of order $\ell$.
\item $\eta:=\det\rho$ is of finite order,
\item $\rho(I_v)$ is finite for $v\in S$, and
\item there exists a place $v\in S$ such that $\rhobar_{v}$ is tame
  and absolutely irreducible and
  $\rho(I_v)\cong\rhobar(I_v)$ --- {\em we say $\rho$ is minimal at $v$}.\label{CondFive} 
\end{enumerate}
Then $R^\eta_X$ is finite flat over~$\BZ_\ell$.
\end{theorem}
The proof of this theorem  will occupy us in the rest of the paper.
Combined with the lifting results of \cite{boecklekhare}, the above
theorem will easily imply Theorem~\ref{SpecialDeJConj} as we see at
the very end (see proof of Theorem \ref{OnDeJConj2}).

Note that if $\ad$ is (absolutely) irreducible, then so is
$\rhobar$. So when we require both assumptions, we only state the
former. 

\bigskip

Let us fix the following notation for this chapter:
For a place $v$ of $\widetilde X$ we denote by $K_v$, $A_v$, $\Fm_v$,
$k_v$ and $q_v$ the completion of the function field $K$ of $X$ at $v$, its
ring of integers, the maximal ideal of the latter, its residue field
and the cardinality of the latter, respectively.
Let $\BA=\BA_K$ be the adeles over $K$. 
The symbols $F$, $\CO$, $\BF$, and $\FM$ will denote a finite extension of
$\BQ_\ell$, its ring of integers, its residue field and the maximal
ideal of $\CO$, respectively. In the definitions to come, $\Lambda$
stands for either of the three rings $F$, $\CO$ or~$\BF$.

Note again that we have chosen a place $w\in X$ and a splitting $s_w$
as in Lemma~\ref{lOnPlaceW} and in the short exact
sequence~(\ref{SESforGabW}). To this splitting corresponds the choice
of a uniformizer $\pi_w$ of $K_w$ unique up to multiplication by
$1+\Fm_w$.

\subsection{Local representations}\label{LocReps}

In this section, we review some background on local representations of
various types. Namely, Galois representations of local Galois groups (or
Weil-Deligne groups), automorphic representations of matrix and
division algebras over a local field, and their
interrelations. Throughout this section we fix a place $v$ of $K$ and
a sufficiently large finite extension $F$ of~$\BQ_\ell$.

\subsubsection{Local automorphic representations}\label{AutoSubSec}

We start by
recalling parts of the classification of absolutely irreducible
$F$-valued representations of the group~$\GL_n(k_v)$. For this let
$k_{v,n}$ be the unique extension of $k_v$ of degree $n$ and
choose a  ring monomorphism $k_{v,n}\into M_n(k_v)$. Let
$\{e_1,\ldots,e_n\}$ be the standard basis of $k_v^n$ and 
$B_n(k_v)\subset\GL_n(k_v)$ the Borel subgroup
which fixes the standard complete flag built from $e_1,\ldots,e_n$.
\begin{definition}\label{ListRepGLn}\begin{enumerate}
\item Denote by $\spec_n$ the Steinberg representation of the group
  $GL_n(k_v)$ defined over $\Lambda$. This is the unique irreducible
  representation of $GL_n(k_v)$ whose dimension is the order of the
  Sylow-$p$ subgroup of $GL_n(k_v)$. 
\item For any character $\chi_v\!:k_{v,n}^*\to \Lambda^*$ such that
  $\chi_v,\chi_v^{q_v},\chi_v^{q_v^2},\ldots, \chi_v^{q_v^{n-1}}$ are
  pairwise distinct, denote by $\theta(\chi_v)$ the representation 
  over $\Lambda$ characterized uniquely by the property
  $$\spec_n\otimes\theta(\chi_v)\cong
  \Ind_{k_{v,n}^*}^{\GL_n(k_v)}\chi_v.$$
\end{enumerate}
We also use the same notation for the representations of the compact
group $\GL_n(A_v)$ obtained by inflation along
$\GL_n(A_v)\onto\GL_n(k_v)$. 
\end{definition}
Note that either of the above representations is defined over $\CO$
if it is defined over $F$. In particular the representation may be
reduced modulo $\Fm$. Moreover if $\chi_v$ is of order prime
to $\ell$, the $\chi_v\pmod \FM$ will have the same order as
$\chi_v$ since reduction modulo $\FM$ on roots of unity of order prime
to $\ell$ is injective. Hence in this case $\theta(\chi_v)\pmod
\FM\cong\theta(\chi_v\pmod\FM)$ is irreducible again.

\begin{lemma} The representations $\spec_n$ and $\theta(\chi_v)$ are
absolutely irreducible. Moreover $\theta(\chi_v)\cong\theta(\chi'_v)$ if
and only if $\chi'_v=\chi_v^{q_v^s}$ for some $s\in\{0,\ldots,n-1\}$.
\end{lemma}

Let $\overline U_0(k_v)\subset \GL_n(k_v)$ be the maximal
parabolic which fixes the subspace $e_1k_v$ and define
$U_0(v):=\{g\in\GL_n(A_v): g\pmod{\Fm_v}\in\overline U_0(v)\}$. 
For $\bar g\in \overline U_0(k_v)$ denote by $a_{11}(\bar g)$ its $(1,1)$-entry. 
\begin{definition}\label{ListRepUz}
For a character $\chi_v\!:k_v^*\to \Lambda^*$, we denote by $I_1(\chi_v)$
the $1$-dimensional representation of $U_0(v)$ defined by
$$U_0(v)\stackrel{\pmod{\Fm_v}}\longto\overline
U_0(k_v)\stackrel{a_{11}}\longto k_v^*
\stackrel{\chi_v}{\longto} \Lambda^*.$$

We define $U_{1,m}(v):=\{g\in U_0(v): a_{11}(g\pmod\FM)\in
k_v^{*\ell^m}\}$, so that $U_{1,m}(v)\subset \kernel(I_1(\chi_v))
\subset U_0(v)$ for any character $\chi_v$ of order~$\ell^m$.
\end{definition}
For later use, we also define $U(v):=1+\Fm_vM_n(A_v)$ and, at the
special place $w$, the compact open group $U_d(w)$ as the kernel of
the composite
$$GL_n(A_w)\stackrel{\mathrm{mod}\,\Fm_w}\onto
\GL_n(k_w)\stackrel{\det}\onto k_w^*\onto k_w^*/(k_w^{*\ell}).$$

\bigskip

Let now $\CD_v$ be a division algebra of degree $n^2$ over its center
$K_v$. Denote by $\CO_{\CD_v}$ its maximal order and by $\Fm_{{\CD_v}}$ the
maximal two-sided ideal of the latter. 
Then $\CO_{{\CD_v}}/\Fm_{{\CD_v}}\cong k_{v,n}$,
so that $\CO_{{\CD_v}}^*/(1+\Fm_{{\CD_v}})\cong k_{v,n}^*$. Thus any
$\chi$ as in Definition~\ref{ListRepGLn}(ii) defines a character
$\widetilde\theta(\chi)$  of finite order of $\CO_{{\CD_v}}^*$ via 
$$\CO_{{\CD_v}}^*\onto k_{v,n}^*\stackrel\chi\longto \Lambda^*.$$
The maximal unramified extension $K_{v,n}$ of $K_v$ of degree $n$ is
contained in ${\CD_v}$, but not uniquely. Therefore the identification 
$\CO_{{\CD_v}}/\Fm_{{\CD_v}}\cong k_{v,n}$ is not unique either and in
fact it is only unique up to an automorphism in
$\Gal(k_{v,n}/k_v)$. One can show:
\begin{prop}\label{ThetaTildeIsom}
There is a bijection between 
\begin{itemize}
\item $\Gal(k_{v,n}/k_v)$-conjugacy
classes of characters $\chi_v$ as in Definition~\ref{ListRepGLn}(ii) and
\item
isomorphism classes of characters
$\widetilde\theta\!:\CO_{{\CD_v}}^*\to \Lambda^*$ up to isomorphisms of~$\CD_v/K_v$
such that $\widetilde\theta$ is trivial on $1+\Fm_{\CD_v}$ and such that
$\widetilde\theta$, $\widetilde\theta^{q_v}$, \ldots,
$\widetilde\theta^{q_v^{n-1}}$ are distinct.
\end{itemize}
\end{prop}

\begin{prop}\label{LocalBadulescu}
Let $\chi_v\!:k_{v,n}^*\to F^*$ be as in
Definition~\ref{ListRepGLn}(ii). Let $\Pi_v$ and $\widetilde\Pi_v$ be smooth
irreducible admissible representations of $\GL_n(K_v)$ and ${\CD_v}^*$,
respectively, so that $\Pi_v$ is square integrable. 
If $\Pi_v^{U(v)}\cong \theta(\chi_v)$ and
$\widetilde\Pi_v^{1+\Fm_{{\CD_v}}}\cong \widetilde\theta(\chi_v)$ as representations
of $\GL_n(A_v)$ and $\CO^*_{{\CD_v}}$, respectively, then under the
local Jacquet-Langlands correspondence (cf.~\cite{badulescux}), one has 
$$\Pi_v\longleftrightarrow\widetilde\Pi_v$$ up to twisting 
by unramified quasicharacters.
\end{prop}
Note that the isomorphism types of both representations do not change
if one replaces $\chi_v$ by~$\chi_v^{q_v}$.

\subsubsection{Local Galois representations}

In this section, we give an explicit description of local residual
representations which are irreducible and at most tamely ramified. 

As before, let $K_{v,n}$ denote the unique unramified extension of $K_v$ of
degree $n$, say with ring of integers $A_{v,n}$. 
Then we have a natural inclusion $I_{v}\into G_{v,n}$
where the latter denotes the absolute Galois group of $K_{v,n}$.
Let $K_{v,n}^*\to G_{v,n}^\ab$ denote the reciprocity map from class
field theory. Its inverse restricted to $I_v$ yields a  morphism 
$\eta'_{v,n}\!: I_v\to A_{v,n}^*$. The composite of this map with the
reduction map $A_{v,n}^*\to k^*_{v,n}$ is defined to be
$$\eta_{v,n}\!: I_v\to k_{v,n}^*.$$ 
\begin{lemma}\label{LocalShapeOfSpecialRho_v}
Suppose $\rhobar_v\!:G_v\to\GL_n(\BF)$ is absolutely irreducible and
at most tamely ramified and assume that all eigenvalues of elements in
$\rhobar_v(I_v)$ lie in $\BF$. Let $(\CC_v,L_v)$ be as in
Proposition 2.16. Then 
\begin{itemize}
\item[(i)] $\rhobar_v$ is of order prime to $\ell$ and 
\item[(ii)] there is a character $\chi_v\!:k_{v,n}^*\to W(\BF)^*$ of
  order prime to $\ell$, unique
  up to the action of $\Gal(K_{v,n}/K_v)\cong\Gal(k_{v,n}/k_v)$, such
  that for any $R\in\CA$ and $\rho_v\in \CC_v(R)$ under the reduction
  map $R\onto \BF$ one has an isomorphism
  $$\left.\rho_v\right|_{I_v}\cong \bigoplus_{i=1}^n
  \chi_v^{q_v^i}\circ\eta_{v,n}.$$
\end{itemize}
\end{lemma}
\begin{remark}
The interested reader may easily prove the following generalization of
(b), which is not needed in the sequel:
Let $\ell^d$ denote the order of an $\ell$-Sylow subgroup of
$k_{v,n}^*$ and set $\CO_j:=W(\BF)[\zeta_{\ell^j}]$ for $0\leq j\leq
d$. Then  for any domain $R\in\CA$ and lift $\rho_v$ of $\rhobar_v$ to
$R$ there exists
$j\in\{0,\ldots,d\}$, a character $\chi_v\!:k_{v,n}^*\to \CO_j^*$, unique up
to the action of $\Gal(k_{v,n}/k_v)$, and an isomorphism
  $$\left.\rho_v\right|_{I_v}\cong \bigoplus_{i=1}^n
  \chi_v^{q_v^i}\circ\eta_{v,n}.$$
\end{remark}
\begin{Proof}
In the present situation, \cite{boecklekhare},~Lem.~2.20, gives
an explicit description of $\rhobar_v$. In the notation of
\cite{boecklekhare},~Lem.~2.20, the irreducibility of
$\rhobar_v$ implies that $d=1$ and that the $U_i$ are matrices of size
$1\times 1$. This proves (i).

Under our hypothesis, \cite{boecklekhare},~Lem.~2.20, moreover
implies that there is a generator of $I_v/P_v$ whose image under
$\rhobar_v$ is a diagonal matrix with distinct diagonal entries
$\bar\mu,\bar\mu^{q_v},\ldots,\bar\mu^{q_v^{n-1}}$, for some
$\bar\mu\in {\BF}$. The description given in the lemma also shows
that $\rhobar_v$ is induced from a $1$-dimensional representation of
$G_{v,n}$, and hence that the action of $I_v$ factors via~$\eta_{v,n}$.

By \cite{boecklekhare},~Rem.~2.18, any lift $\rho_v$ in $\CC_v(R)$, $R\in\CA$,
satisfies $\rho_v(I_v)\stackrel\cong\to\rhobar_v(I_v)$. Therefore the image of
the above generator of $I_v/P_v$ under $\rho_v$ is still
diagonal of the same order. By functoriality, it suffices to consider
the case $R=W(\BF)$, in which (ii) follows from a suitable choice
of~$\chi_v$.
\end{Proof}

For $\rho_v$ as in the previous lemma, the group $\rho_v(I_v)$ is
finite, and so $\rho_v$ agrees with its Deligne-Weil representation.
We have the following standard result that is part of the local Langlands correspondence
for function fields proved in \cite{lrs}:
\begin{prop}\label{LocalLanglForChi}
Suppose $\rhobar_v$ is as in the previous lemma, $\rho_v\!:G_v\to
\GL_n(\CO)$ is a deformation of type $\CC_v$, and $\chi_v$ lies in the
corresponding conjugacy class of  characters. Let $\Pi_v$ denote an
irreducible smooth admissible representation of $\GL_n(K_v)$.
Then $\rho_v \leftrightarrow \Pi_v$ under the local Langlands
correspondence, up to twisting 
by unramified quasicharacters, if and only if there is an isomorphism $\Pi_v^{U(v)}
\cong \theta(\chi_v)$ as a representation of~$\GL_n(A_v)$.
\end{prop}

\subsection{Automorphic forms over function fields}\label{AutoBasics}

Let us fix a division algebra $\CD$ which is central over $K$ of
degree $n^2$. Our next aim is
to define certain spaces of cusp forms on $\GL_n(\BA)$ and $\CD^*(\BA)$.

\medskip

In the sequel we will often work under the following hypothesis:
\begin{assumption}\label{Aut-BasicGlAss} 
\begin{itemize}
\item[(a)] $\# S\geq2$,
\item[(b)] $\CD$ is totally split at all places not in $S$ and (maximally)
  ramified of degree $n$ at the places in $S$ (such a $\CD$ exists by class field theory),
\item[(c)] $T$ is a finite subset of $X\setminus\{w\}$,
\item[(d)] $\omega\!:\GL_1(\BA)\to \Lambda^*$ is a character of
  finite order, trivial on~$K^*$, unramified outside $S$, and at most
  tamely ramified,
\item[(e)] for each $v\in S$, $\chi_v\!:k_{v,n}^*\to \Lambda^*$ is a
  character such that the characters $\chi^{q_v^i}_v$,
  $i=1,\ldots,n$, are pairwise distinct,
\item[(f)] for $v\in S$, the characters $\chi^{1+q_v+\ldots+q_v^{n-1}}\!:k_v^*\cong
  A_v^*/(1+\Fm_v)\to \Lambda^*$ and $\left.\omega\right|_{A_v^*}$
  agree as characters of $A_v^*$, and
\item[(g)] for $v\in T$, $\chi_v\!:k_v^*\to \Lambda^*$ is a character 
  of $\ell$-power order $\ell^{m_v}$, which may be trivial
\end{itemize}
\end{assumption}
The above set-up will be of use in two different instances, namely in
lowering the level of an automorphic cusp form associated to a
residual Galois representation and in constructing Taylor-Wiles
systems. In the first case, we will use the above notation as stated,
in the second, the set $T$ will be denoted by $Q$ or $Q_m$,
respectively.

\smallskip

Let $\Um$ denote the tuple of the $m_v$, $v\in T$, and let
$m\in\BN_0$ satisfy $m\leq\min\{m_v:v\in T\}$. Using the chosen
uniformizer $\pi_w\in K_w$ corresponding to $s_w$, we define
$$Z^{\Um}_{T}:=K^*\cdot\Big(\pi_w^\BZ(1+\Fm_w)k_w^{*\ell}\times
\prod_{v\in T}(1+\Fm_v)k_v^{*\ell^{m_v}}\times\prod_{v\notin
  T\cup\{w\}} A_v^* \Big)\subset\GL_1(\BA)$$
and $Z^m_{T}:=K^*(\pi_w^\BZ(1+\Fm_w)k_w^{*\ell}\times
\prod_{v\in T}(1+\Fm_v)k_v^{*\ell^m}\times\prod_{v\notin
  T\cup\{w\}} A_v^*)$.

\smallskip

Under
the above hypothesis we define compact open subgroups of $\GL_n(\BA)$
and $\CD^*(\BA)$ by 
\begin{eqnarray*}
U_T^\Um&:=& \prod_{v\in S} U(v) \times U_d(w)\times \prod_{v\in T}
U_{1,m_v}(v) \times\!\!\!\!\prod_{v\notin S\cup T\cup\{w\}}\!\!\!\GL_n(A_v)\subset \GL_n(\BA),\\ 
\widetilde U^\Um_T&:=& \prod_{v\in S} (1+\Fm_{\CD_v})\times U_d(w)\times\prod_{v\in
  T} U_{1,m_v}(v) \times\!\!\!\!\prod_{v\notin S\cup T\cup\{w\}}\!\!\!\GL_n(A_v)\subset\CD^*(\BA),\\ 
\widetilde U_T^m&:=& \prod_{v\in S} (1+\Fm_{\CD_v})\times U_d(w)\times\prod_{v\in
  T} U_{1,m}(v) \times\!\!\!\!\prod_{v\notin S\cup T\cup\{w\}}\!\!\!\GL_n(A_v)\subset\CD^*(\BA),\\
\widetilde U^{00}_T&:=& \prod_{v\in S} \CO_{\CD_v}^*\times U_d(w)\times\prod_{v\in
  T} U_0(v) \times\!\!\!\!\prod_{v\notin S\cup T\cup\{w\}}\!\!\!\GL_n(A_v)\subset\CD^*(\BA).
\end{eqnarray*}

\smallskip

We define $\pbC^{\omega,(\chi_v)}_{S,T}(\Lambda)$ as the space
of all functions 
$$f\!:\GL_n(K)\backslash \GL_n(\BA)/U_T^\Um \to \Lambda
$$ 
with the following properties:
\begin{enumerate}
\item The central action of $Z_{T}^\Um$ on $f$ is described by $\omega$,
\item for $v\in S$, the right action of $\GL_n(A_v)$ on the
  finite-dimensional free $\Lambda$-module
  $f(g\cdot\ublank)\Lambda[\GL_n(A_v)]$ is via $\theta(\chi_v)$,
\item for $v\in T$ the right action of $U_{1,m_v}(v)$ on
  $f(g\cdot\ublank)$ is via the character $I_1(\chi_v)$, and
  \label{CondthreeOnCuspForms} 
\item $f$ is cuspidal (cf. \cite{boreljacquet}, \S~5).
\end{enumerate}
Note that while condition (i) makes for $v\in S\cup T$ an assertion on
the central action of $(1+\Fm_v)k_v^{*\ell^{?}}$ only for some
$?\in\BN$, conditions (ii) and (iii) do completely determine the
central action of~$A_v^*$.

\medskip

Similarly, we define the space $\OpbC^{\omega,(\chi_v)}_{S,T}(\Lambda)$
as the space of all functions  
$$f\!:\CD^*(K)\backslash \CD^*(\BA)/\widetilde U_T^\Um \to \Lambda$$ 
with the following properties:
\begin{enumerate}
\item For $z\in Z_{T}^\Um$ one has $f(zg)=\omega(z)f(g)$,
\item for $v\in S$ and $h\in\CO^*_{\CD_v}$ one has $f(gh)=f(g)\tilde\theta(\chi_v)(h)$,
\item for $v\in T$ and $u\in U_{1,m_v}(v)$ one has
  $f(gu)=f(g)I_1(\chi_v)(u)$.
\end{enumerate}
Since $\CD^*(\BA)$ is compact modulo its center, we do not need to
impose any cuspidality conditions on
$\OpbC^{\omega,(\chi_v)}_{S,T}(\Lambda)$.

For $\CD$ we define a second space of functions
$\OpbC^{\omega,m}_{S,T}(\Lambda)$ as 
the space of all functions  
$$f\!:\CD^*(K)\backslash \CD^*(\BA)/\widetilde U_T^m \to \Lambda$$ 
with the following properties:
\begin{enumerate}
\item For $z\in Z_{T}^m$ one has $f(zg)=\omega(z)f(g)$,
\item for $v\in S$ and $h\in\GL_n(A_v)$ one has $f(gh)=f(g)\theta(\chi_v)(h)$,
\end{enumerate}

For any space $\OpbC^{*}_{?}(\Lambda)$, we define a
corresponding space of cusp forms $\ObC^{*}_{?}(\Lambda)$ as the
quotient of the former space by its subspace of functions that factor
through the norm. 
The spaces $\OpbC^{*}_{?}(\Lambda)$ play mainly an
intermediary role to prove results about the spaces
$\ObC^{*}_{?}(\Lambda)$ which are the ones we are interested in.

The following is an immediate consequence of the above definitions:
\begin{prop}\label{Aut-DirSumOfCuspForms}
Let $m\in\BN$ be such that $\ell^m$ divides the order of
$k_v^*$ for all $v\in T$ and assume that $F$ contains
$W(\BF)[\zeta_{\ell^m}]$. Then 
$$\OpbC^{\omega,m}_{S,T}(F)\cong\bigoplus_{(\chi'_v)}
\OpbC^{\omega,(\chi'_v)}_{S,T}(F),$$ 
where the $(\chi'_v)$ range over all sets of characters such that
$\chi_v'=\chi_v$ for $v\in S$ and such that the order of $\chi_v'$
divides $\ell^m$ for $v\in T$.
\end{prop}

\begin{prop}\label{PropRedOfCuspForms}
Let $M(\CO)$ be any of the spaces
$\pbC^{\omega,(\chi_v)}_{S,T}(\CO)$,
$\OpbC^{\omega,(\chi_v)}_{S,T}(\CO)$ or
$\OpbC^{\omega,m}_{S,T}(\CO)$.
\begin{enumerate}
\item \label{PROCFi} The spaces $M(\CO)$ are free and finitely
  generated over~$\CO$. 
\item \label{PROCFii} The induced morphism $ M(\CO)\otimes_\CO\BF
  \to M(\BF)$ is injective.
\end{enumerate}
\end{prop}
We suspect that the morphism in \ref{PROCFii} is not always
surjective. But we do expect surjectivity after 
localizing at a maximal ideal of the Hecke algebra (as in the next
section) for which the corresponding $n$-dimensional mod $\ell$ Galois
representation is irreducible. That we cannot prove this is the main
reason why we have local conditions on $\rhobar$ in the main theorem
of this work. 
 
\begin{Proof}
Clearly the spaces $M(\CO)$ inject into $M(F)$ and since the latter
are finite-dimensional, \ref{PROCFi} is immediate.

For \ref{PROCFii} note that cuspidality is preserved under the
reduction map $\CO\to\BF$ and that cusp forms are functions on a set.
\end{Proof}

\medskip

By their very definition, cf.\ \cite{boreljacquet}, one has a smooth
admissible automorphic representation $\Pi(f)$ of $\GL_n(\BA)$ attached to
an automorphic form $f$ for $\GL_n$ (simply given by $\Lambda[\GL_n(\BA)]f$). 
The local constituents of $\Pi(f)$ are denoted $\Pi_v(f)$, so that
$\Pi(f)\cong\hat\otimes_v \Pi_v(f)$. Conversely, if $\Pi$ is a smooth
admissible cuspidal automorphic representation and $U$ a compact open
subgroup of $\GL_n(\BA)$, then $\Pi^U$ is a (possibly empty) space of cusp
forms. 

Viewing $\CD^*$ as an algebraic group over $K$, again following
\cite{boreljacquet}, we denote by $\widetilde\Pi(\tilde f)$ the smooth
admissible automorphic representation of $\CD^*(\BA)$ attached to an
automorphic form $\tilde f$ for $\CD^*$. We write 
$\widetilde\Pi_v(\tilde f)$ for the local constituents at $v$. Again
one has $\widetilde\Pi(\tilde f)\cong\hat\otimes_v\widetilde\Pi_v(\tilde f)$
and forms arise as fixed spaces of smooth admissible irreducible
representations under compact open subgroups. 

\medskip

We state a special case of the global Jacquet-Langlands
correspondence in the function field case, as worked out by Badulescu,
\cite{badulescuy}, Prop.~4.7.
\begin{theorem}\label{GlobalThmBadu}
Suppose the hypothesis of Assumption~\ref{Aut-BasicGlAss} hold for
$\Lambda$ a field $F$. Then to every irreducible cuspidal automorphic
representation $\Pi=\hat\otimes_v\Pi_v$ of $\GL_n(\BA)$ such that
$\Pi_v^{U(v)}\cong \theta(\chi_v)$ for $v\in S$, there exists a unique
irreducible cuspidal automorphic representation
$\widetilde\Pi=\hat\otimes_v\widetilde\Pi_v$ of $\CD^*(\BA)$ such that
$\widetilde\Pi_v^{1+\Fm_{\CO_{\CD_v}}}\cong \widetilde\theta(\chi_v)$ for
$v\in S$ and $\Pi_v\cong\widetilde\Pi_v$ for all $v\notin S$, and vice
versa. The correspondence is compatible with the local
Jacquet-Langlands correspondence, and thereby uniquely determined. 
\end{theorem}

\subsection{Hecke-algebras}\label{Hecke}

An important tool in the theory of automorphic forms is the action of
suitable Hecke-algebras. For $v\notin S\cup T\cup\{w\}$ we define
$\CH_v$ as the Hecke-algebra of bi-$\GL_n(A_v)$-invariant locally
constant compactly supported functions on $\GL_n(K_v)$ with
values in $\BZ_\ell$, and where multiplication is given by
convolution. The algebra $\CH_v$ contains naturally defined elements
$T_{v,i}$, $i=1,\ldots,n$, the Hecke-operators at $v$, and the Satake
isomorphism asserts that
$$\CH_v\cong \BZ_\ell[T_{v,1},\ldots,T_{v,n},T_{v,n}^{-1}].$$
Via integration, one defines a compatible action of $\CH_v$ on
automorphic forms and representations. The resulting action is only on
the $v$-th component of the the adelic group $G(\BA)$,
$G\in\{\GL_n,\CD^*\}$. It preserves the spaces of cusp forms defined
above.

\smallskip

For $v\in T$ ($T$ might be empty) and $m_v\in\BN_0$, we follow \cite{harristaylor}, Ch.~II
(there is no difference between the function and the number field case
here): So we denote by $V_{v,i}$, $i=1,\ldots,n-1$, the Hecke-operators
in the convolution algebra of $U_{1,m_v}(v)$-bi-invariant locally
constant compactly supported functions on $\GL_n(K_v)$ with
values in $\BZ_\ell$, defined as in \cite{harristaylor},
(II.2.2). By $V_{v,n}$ we denote the {\em $U$-operator} from
\cite{harristaylor}, II.2.4, which again lies in the above convolution
algebra. (We choose the notation $V_{v,n}$, to avoid any conflict with
our notation for compact opens.) The commutative subalgebra spanned by
the $V_{v,i}$, $i=1,\ldots,n$, is denoted by~$\CH^{m_v}_{v}$.  We
consider Hecke action at places in $T$ only when building Taylor-Wiles
systems.

\medskip

To obtain a Hecke action at almost all places (namely outside $S\cup\{w\}$), we set
$$\CH_T^\Um:=
\bigotimes_{v\in T}\CH^{m_v}_v\otimes
\bigotimes_{v\notin S\cup\{w\}\cup T}\CH_v.$$ 
Again the operation of $\CH_T^\Um$ preserves the spaces of cusp forms
$\pbC^{\omega,(\chi_v)}_{S,T}(\Lambda)$,
$\ObC^{\omega,*}_{S,T}(\Lambda)$  and
$\OpbC^{\omega,*}_{S,T}(\Lambda)$. For $0\leq m\leq\min\{m_v:v\in
T\}$, we define
$\widetilde\CH_T^{m}(\Lambda)$ as the image of
$\CH_T^\Um\otimes_\BZ\Lambda$ in the endomorphism ring of
$\ObC^{\omega,m}_{S,T}(\Lambda)$. 
Because $\ObC^{\omega,m}_{S,T}(\Lambda)$ is a free $\Lambda$-module of
finite rank, the same holds for $\widetilde\CH_T^{m}(\Lambda)$ and
one has
$\widetilde\CH_T^{m}(F)\cong\widetilde\CH_T^{m}(\CO)\otimes_\CO F$.

\medskip

We state the main theorem of \cite{lafforgue}:
\begin{theorem}[\cite{lafforgue}]\label{LaffThm}
For any finite subset $T$ of $X$, there is a bijection between 
\begin{itemize}
\item smooth irreducible cuspidal automorphic representations $\Pi$
  whose central character is of finite order and with $\Pi_v$
  unramified for $v\notin S\cup\{w\}\cup T$, and  
\item irreducible continuous representations
  $\rho\!:\pi_1(X\setminus (\{w\}\cup T))\to\GL_n(\overline \BQ_\ell)$ with determinant
  of finite order.
\end{itemize}
Suppose that for $\Pi$ as above, the eigenvalues of the operators
$T_{v,i}$, $v\notin S$, $i=1,\ldots,n$, are $\alpha_{v,i}$ (and
$\alpha_{v,0}=1$). Then the correspondence is given by the condition
\begin{equation}\label{GLcorrespEq}
\det\Big(1-x \rho(\Frob_v)\Big)=\sum_{i=0}^n x^{n-i}\alpha_{v,i}.\end{equation}
\end{theorem}

Suppose we are given $\rho\!:\pi_1(X\setminus \{w\}\cup
T)\to\GL_n(\CO)$ with absolutely irreducible residual representation
$\rhobar$. Let $\Pi$ be the corresponding
automorphic representation. 
By the relation (\ref{GLcorrespEq}) the Hecke-eigenvalues
$\alpha_{v,i}$, $v\notin S\cup\{w\}\cup T$, $i\in 1,\ldots,n$, 
lie in $\CO$, and so they define ring-homomorphisms
$$\CH_v\to\BF:T_{v,i}\mapsto\ \alpha_{v,i}\!\!\!\!\pmod \FM$$
for $v\notin S\cup\{w\}\cup T$. 

If for $v\in T$ the representation $\Pi_v$ satisfies
$\Pi_v^{U_{1,m_v}}\neq0$, if $\rhobar$ is unramified at $v$ with
$n$ distinct eigenvalues and if $\bar\lambda$ is the eigenvalue at the
$(1,1)$-entry, then by \cite{harristaylor} one has the following: $\Pi_v$
is principal series with ramification at most at the character
$\chi_v$ of the $(1,1)$-entry, and of order dividing $\ell^{m_v}$.
There are now two cases:
\label{HTquote} if $\Pi_v$ is unramified, then
$\dim\Pi_v^{U_{1,m_v}}=n$, the action of $V_{v,n}$ on this space
is semisimple, it has a unique eigenvalue which reduces to
$\bar\lambda$ and the corresponding eigenspaces is
$1$-dimensional; this $1$-dimensional subspace is also an eigenspace
for the operators $V_{v,i}$, $i=1,\ldots,n-1$. The corresponding
eigenvalues are denoted $b_{v,i}$ and their mod $\FM$ reductions
depend only on~$\rhobar(\Frob_v)$.

If on the other hand, $\Pi_v$ is ramified, then
$\dim\Pi_v^{U_{1,m_v}}=1$, and if, we denote by $b_{v,i}$ the
eigenvalues for the Hecke-operators $V_{v,i}$, $i=1,\ldots,n$, their
reductions modulo $\FM$ depend only on $\rhobar(\Frob_v)$ (and by the
same formulas as in the first case). Moreover the reduction of
$b_{v,n}$ is $\bar\lambda$. So in either case one obtains a ring
homomorphism
$$\CH_v^{m_v}\to \BF: V_{v,n}\mapsto\bar\lambda, V_{v,i}\mapsto
b_{v,i}\!\!\!\pmod\FM \,\,\,(i=1,\ldots,n-1).$$

\smallskip

The above homomorphisms induce a ring homomorphism $\CH_T^{\Um}\to
\BF$ whose kernel is a maximal ideal which we denote by
$\FM_{\rhobar}$. This notation is justified because the homomorphism
depends only on data defined in terms of~$\rhobar$. We also denote by
$\FM_{\rhobar}$ the image of this ideal in 
$\widetilde{\CH}_T^m(\CO)$. (Since a priori there is no relation
between $\rhobar$ and $\ObC_{S,T}^{\omega,m}(\CO)$, this
image can be all of $\widetilde\CH_T^{m}(\CO)$.) For an
$\CH_T^\Um$-module $M$, we denote by $M_{\FM_{\rhobar}}$ its localization
at~$\FM_{\rhobar}$. 

The following lemma is immediate from the fact that $\rhobar$ is irreducible
and the fact that on functions in either $\OpbC^{\omega,(\chi_v)}_{S,T}(\CO)$ or
$\OpbC^{\omega,m}_{S,T}(\CO)$ that factor through the norm, Hecke
operators acts by just scalars that are their degrees.
\begin{lemma}\label{stupid} 
We have isomorphisms 
$\OpbC^{\omega,(\chi_v)}_{S,T}(\CO)_{\FM_{\rhobar}} \simeq \ObC^{\omega,(\chi_v)}_{S,T}(\CO)_{\FM_{\rhobar}}$ and 
$\OpbC^{\omega,m}_{S,T}(\CO)_{\FM_{\rhobar}} \simeq \ObC^{\omega,m}_{S,T}(\CO)_{\FM_{\rhobar}}$
compatible with the Hecke action of~$\CH_T^\Um$.
\end{lemma}

Since for $v\in  T$, the global Jacquet-Langlands correspondence
given in Theorem~\ref{GlobalThmBadu} is locally the identity, we
have the following corollary:
\begin{cor}
The Jacquet-Langlands correspondence in Theorem~\ref{GlobalThmBadu} is
compatible with the Hecke-action of~$\CH_T^\Um$.
\end{cor}

This has the following 2 immediate corollaries:
\begin{cor}\label{FirstCorOnBadu}
Under Assumption~\ref{Aut-BasicGlAss} there is a natural  isomorphism
of $\CH_T^\Um$-modules 
$$\Hom_{\prod_{v\in S}\GL_n(k_v)}\Big(\bC^{\omega,(\chi_v)}_{S,T}(F),
\prod_{v\in S} \theta(\chi_v)\Big)\stackrel\cong\longto
\ObC^{\omega,(\chi_v)}_{S,T}(F).$$
\end{cor}
\begin{Proof}
There is a bijection between irreducible constituents $\tau$ of
$\bC^{\omega,(\chi_v)}_{S,T}(F)$  as a representation of
$\prod_{v\in S}\GL_n(k_v)\times \prod_{v\in T} \overline
U_0(v)$ and irreducible admissible representations $\Pi$ such that
$\Pi\cong\Pi(f)$ for some $f\in
\bC^{\omega,(\chi_v)}_{S,T}(F)$. Both types of representations
are preserved under the Hecke-action of $\CH_T^\Um$, and the
correspondence preserves Hecke-eigenvalues. Moreover as a
representation of $\GL_n(\kappa_v)$, $v\in S$, the isomorphism type of
any $\tau$ is $\theta(\chi_v)$. For a place $v$ not in $S$, the
representation type of $\tau$ is a $1$-dimensional representation. In
particular, the dimension of 
$$V_1:=\Hom_{\prod_{v\in S}\GL_n(\kappa_v)}\Big(
\bC^{\omega,(\chi_v)}_{S,T}(F), \prod_{v\in S}
\theta(\chi_v)\Big)$$
is the number of isomorphism classes of representation $\Pi$. Moreover
any system of Hecke-eigenvalues of such a $\Pi$ occurs as the
Hecke-eigensystem of a $1$-dimensional Hecke eigen-subspace of~$V_1$.

A similar reasoning yields the analogous assertion for 
$$V_2:=\Hom_{\prod_{v\in S}\GL_n(\kappa_v)}\Big(
\ObC^{\omega,(\chi_v)}_{S,T}(F),\prod_{v\in S}
\widetilde\theta(\chi_v)\Big).$$
By the above theorem one has $V_1\cong V_2$ as
Hecke-modules. Because each of the $\widetilde\theta(\chi_v)$ is
$1$-dimensional, the assertion follows.
\end{Proof}
\begin{cor}\label{SecondCorOnBadu}
Under Assumption~\ref{Aut-BasicGlAss} there is a natural  isomorphism
of $\CH_T^\Um$-modules 
$$\Hom_{\prod_{v\in S}\GL_n(\kappa_v)}\Big(
\bC^{\omega,(\chi_v)}_{S,T}(\CO)_{\FM_{\rhobar}},
\prod_{v\in S} \theta(\chi_v) \Big)\stackrel\cong\longto
\ObC^{\omega,(\chi_v)}_{S,T}(\CO)_{\FM_{\rhobar}}\otimes_\CO F.$$
\end{cor}
Note that the representations $\theta(\chi_v)$ are $F$-valued.
\begin{Proof}
All Hecke-eigenvalues are integral over $\CO$. Therefore the
Hecke-modules in Corollary~\ref{SecondCorOnBadu} are the direct sum of
those Hecke-eigenspaces of the modules in
Corollary~\ref{FirstCorOnBadu} (with
$F$-coefficients) for which the
$\CO$-integral eigenvalues satisfy the congruence conditions
described by $\FM_{\rhobar}$. 
\end{Proof}\label{ThirdCorOnBadu}

\subsection{Carayol's principle}\label{carprin}

The methods of this section, which allow us to change ``types'' of automorphic forms 
which give rise to a given $\rhobar$, use a principle discovered by
Carayol that occurs in the proof of \cite{carayol}~Lemme~1.

We define $N_n(\ell)$ as the order of an $\ell$-Sylow subgroup of
$\GL_n(k)$. 

{\it We may consider the Hecke action without the operators at places in $T$ and will denote
the induced maximal ideal of this smaller Hecke algebra by the same symbol.}

\begin{lemma}\label{AutoLoweringLevel}
We assume the set-up of Assumption~\ref{Aut-BasicGlAss} for sets of
characters $(\chi_v)$ and $(\chi_v')$ such that
\begin{itemize}
\item[(a)] $\Lambda$ is the fraction field $F$ of $\CO$,
\item[(b)] $\chi_v=\chi_v'$ for $v\in S$,
\item[(c)] for each $v\in T$, the product of $N_n(\ell)$ with the order
  of $\chi_v^{-1}\chi_v'$, which is a power of $\ell$, divides the
  order of~$k_v^*$.
\end{itemize}
Then 
$$\rank_\CO \ObC^{\omega,(\chi_v)}_{S,T}({\cal O})_{\FM_{\rhobar}}=\rank_\CO
\ObC^{\omega,(\chi'_v)}_{S,T}({\cal O})_{\FM_{\rhobar}}$$ (and the same conclusion holds if we consider
the corresponding maximal ideal of the smaller Hecke algebra without the operators at $T$).
\end{lemma}

Localization at $\FM_{\rhobar}$ commutes with reduction $\CO\to\BF$ of
the ring of integers of $F$ to its residue field. Furthermore, after
reduction on has $\chi_v\equiv\chi_v'\equiv1$ in $\BF$ for all $v\in
T$. So in view of Proposition~\ref{PropRedOfCuspForms}(i)-(ii) and Lemma \ref{stupid}, it
suffices to prove the following
\begin{lemma}\label{crucial}
Suppose the hypothesis of Lemma~\ref{AutoLoweringLevel} hold. Then the submodules 
$\OpbC^{\omega,(\chi_v)}_{S,T}(\CO)\otimes_{\CO}\BF$ and
$\OpbC^{\omega,(\chi'_v)}_{S,T}(\CO)\otimes_{\CO}\BF$
of $\OpbC^{\omega,(\chi'_v)}_{S,T}(\BF)$ agree. 
\end{lemma}

\begin{Proof}
By $Z$ we denote the center of $\CD^*$, viewed as an algebraic group
over $K$, i.e., $Z\cong\BG_m$, and we regard $Z_T^\Um$ as a subgroup
of $Z(\BA)$. Because of the central action and the
conditions on the places $v\in S\cup T$, elements in
$\OpbC^{\omega,(\chi_v)}_{S,T}(\CO)$ maybe thought of as functions
on the finite set $\CD^*Z_{T}^\Um\backslash\CD^*(\BA)/\widetilde U^{00}_T$.
Therefore we choose finitely many $d_j\in\CD^*(\BA)$, $j\in J$, such that
$$\CD^*(\BA)= \coprod_{j\in J} \CD^*Z_{T}^\Um d_j \widetilde U^{00}_T.$$
Let $J_0$ denote the subset of those $j\in J$ such that $g(d_j)=0$ for all
$g\in \OpbC^{\omega,(\chi_v)}_{S,T}(\CO)$, and define
analogously $J'_0$ for the $\chi'_v$. We claim that
\begin{enumerate}
\item For any choice $\alpha_j\in\CO$, $j\in J\setminus J_0$,
  there is a (unique) $g\in \ObC^{\omega,(\chi_v)}_{S,T}(\CO)$
  with $g(d_j)=\alpha_j$, and the analogous assertion in the primed
  setting.
\item $J_0=J'_0$.
\end{enumerate}
This clearly implies the lemma. 

Let $Z^\Um_{T,0}\subset Z_{T}^\Um$ denote the kernel of the degree map on
$Z_{T}^\Um$, i.e., the subgroup obtained from $Z_T^\Um$ by replacing the
$w$-component by $(1+\Fm_w)k_w^{*\ell^{m_v}}\subset
\pi_w^\BZ(1+\Fm_w)k_w^{*\ell^{m_v}}$. 
In the sequel we prefer to work with the compact group
$Z^\Um_{T,0}/Z(K)$ rather than with $Z^\Um_{T,0}$. Therefore, since any cusp
form is constant along any $Z(K)$-orbit, we will regard
cusp forms as functions on $\CD^*(\BA)/Z(K)$.

Elements in $\CD^*(K)/Z(K)$ are denoted $\gamma,\gamma'$, elements in
$Z_{T}^\Um/Z(K)$ are denoted $z,z'$ and those in $U^{00}_TZ(K)/Z(K)\cong
U^{00}_T/Z(k^*)$ by $u,u'$. By $u_v$ and $u_v'$, respectively, we denote
the components at $v$ in $\CO^*_\CD/k^*$.
Suppose $g\in \ObC^{\omega,(\chi_v)}_{S,T}(\CO)$ and
$\alpha_j=g(d_j)$. Then if $x=\gamma z d_j u =\gamma' z' d_j u'$ we deduce
$$\omega(z)\prod_{v\in S\cup T}\chi_v(u_v)\,\alpha_j=g(x)=
\omega(z')\prod_{v\in S\cup T}\chi_v(u'_v)\,\alpha_j.$$ 
Moving the $\gamma$'s to one side and the other variables to the
other, we deduce the same equality from $({\gamma'}^{-1}\gamma d_j)
=d_j (z^{-1}z' )(u'u^{-1})$. We may apply the degree map after taking
determinants, since the degree map vanishes on $K^*$. This yields
$\degr(z^{-1}z')=\degr({u'}^{-1}u)=0$.

Thus we have $j\in J_0$ if and only if there exists $z\in
Z^\Um_{T,0}/Z(K)$, $\gamma\in \CD^*(K)/Z(K)$ and $u\in U_T^{00}/k^*$ such that 
$$\gamma d_j=d_j z u\hbox{ and }\omega(z)\prod_{v\in S\cup T}\chi_v(u_v)\neq1.$$
This proves (i) above.

\medskip

To approach (ii), we first introduce some notation. For $v\in T$
define $m_v\in\BZ$ by the condition that $\ell^{m_v}$ is the
$\ell$-part of the order of $k_v^*$. By hypothesis (c) we have
$\ell\cdot N_n(\ell)| \ell^{m_v}$, and we define
$\Um$ as the tuple consisting of the $m_v$. We consider the finite group 
$$M^\Um_T:=\Big(Z^\Um_{T,0}\widetilde U_T^{00}/Z(K) \Big)/
\Big(\widetilde U_{T}^\Um Z(K)/Z(K) \Big).$$ 
Note that $\omega$ is a character on $Z(\BA)/Z(K)$, while $\chi_v$
is a character on the $v$-component of $\widetilde U_T^{00}/\widetilde 
U_{T}^\Um$ for $v\in S
\cup T$. By Assumption~\ref{Aut-BasicGlAss}(f), $\omega$
is compatible with the central action on $\theta(\chi_v)$. For $v\in
T$, the characters $\omega$ as well as $I_1(\chi_v)$ if restricted to
elements in $(1+\Fm_v)k_v^{*\ell^{m_v}}$ are trivial. Therefore there 
exists a unique character on $M_T^\Um$, we call it $\chi_T$, which
agrees with $\omega$ and the $\chi_v$, when restricted to the
appropriate subgroups of~$M_T^\Um$.

Define $H_j:=d_j^{-1}\CD^*(K)/Z(K)d_j\cap U_T^{00}Z^\Um_{T,0}/Z(K)$. Since
the left hand side is discrete and the right hand side is compact,
$H_j$ is finite, and we denote by $\overline H_j$ its image in
$M_T^\Um$. We then have
$$J_0=\{j\in J : \chi_T \hbox{ is non-trivial on } \overline H_j\}$$
and a similar characterization of $J_0'$ with the analogously
defined~$\chi'_T$.

By the lemma below, the $\ell$-part of the exponent of $H_j$ and hence
also of $\overline H_j$ divides $N_n(\ell)$. Then  by hypothesis (c) for
any $v\in T$ the image of $\overline H_j$ under 
$$M_T^\Um\onto \widetilde U_0(v)/\widetilde U_{1,m_v}(v)\cong\BZ/\ell^{m_v}\BZ \onto
\BZ/\ell \BZ$$ 
is trivial, and so $\overline H_j$ lies in the kernel of
$\chi_v^{-1}\chi_v'$. This means that $\chi_T$ and $\chi'_T$ agree on
$\overline H_j$ for all $j$, and thus part~(ii) of the above claim is
shown. 
\end{Proof}
\begin{lemma}
If $\gamma\in\CD^*/K^*$ is of $\ell$-power order, its order 
divides $N_n(\ell)$.
\end{lemma}
\begin{Proof}
Let $S\CD^*$ be the kernel of the reduced norm on $\CD^*$. Then the
short exact sequence $1\to S\CD^*\to \CD^*\to K^*\to 1$ induces an
exact sequence
$$0\longto S\CD^*/\{x\in K^*:x^n=1\}\longto \CD^*/K^*\longto
K^*/K^{*n}\longto 0.$$ 
Because the image of $\gamma$ in $K^*/K^{*n}$ is of $\ell$-power order,
it must lie in $k^*/k^{*n}$. Because $\ell$ is prime to $n$, the order of
this image divides the order of the $\ell$-Sylow of $k^*$. Say the latter
order is $n_\ell$. Then $\gamma^{n_\ell}$ lies in $S\CD^*/
\{x\in K^*:x^n=1\}$. Again because $n$ is prime to $\ell$, there is an
element $\gamma'\in S\CD^*$ which maps to $\gamma^{n_\ell}$ and whose
order is the same as that of $\gamma^{n_\ell}$. Thus it suffices to show
that if $\gamma'\in S\CD^*$ is of $\ell$-power order $\ell^d$, then the order of
$\gamma'$ divides the order of an $\ell$-Sylow of $\SL_n(k)$.

Let now $K'/K$ be a splitting field of $\CD$, say of degree $n$, whose
constant field is still $k$. Then we may regard $\gamma'$ as an
element of $\SL_n(K')$. Let $k''$ be the smallest extension of $k$
which contains a primitive $\ell^d$-th root of unity, and let $K''\supset
K'$ be the corresponding unramified extension. We assume that
$\gamma'$ is given in rational canonical form over $K'$. We claim that
$\gamma'$ has entries in $k$. For this we may assume that the rational
canonical form consists of a single block, which is thus completely
determined by the characteristic polynomial of $\gamma'$. But the
characteristic polynomial has coefficients in $k''$ and in $K'$ and
thus in $k$. Hence $\gamma'$ has coefficients in $k$. Thus $\gamma'$
lies up to conjugation in $\SL_n(k)$, and our last assertion is shown.
\end{Proof}

\subsection{Taylor Wiles systems}\label{TWSyst}

Throughout this section, we fix a lift $\rho\!:\pi_1(X)\to \GL_n(\CO)$ of
$\rhobar$ and make the following assumptions:
\begin{assumption}\label{TWsysAss}
\begin{itemize}
\item[(a)] $\eta:=\det\rho$ is of finite order,
\item[(b)] for all $v\in S$, $\rhobar_v$ is tamely ramified and
  absolutely irreducible, and $\rho_v(I_v)\to\rhobar_v(I_v)$ is an
  isomorphism.
\end{itemize}
\end{assumption}

For each $m\in\BN$ we also fix a finite subset $Q_m \subset
X\setminus\{w\}$ of places such that for all $v\in Q_m$ the matrix
$\rhobar(\Frob_v)$ has distinct eigenvalues and $q_v\equiv 1 \pmod
{\ell^m}$. In this section we will complement the Galois theoretic work
in Section 2 by automorphic results that together yield the existence
of Taylor-Wiles systems.

Define $\Delta_m$ as the maximal quotient of
$\prod_{v\in Q_m}k_v^*$ of exponent $\ell^m$. Via the projection onto
the $(1,1)$-entry, we identify 
$$\Delta_m \cong \prod_{v\in Q_m} U_0(v)/U_{1,m}(v).$$
The latter group naturally acts on $\widetilde U_{Q_m}^0/\widetilde U_{Q_m}^m$ (this is the ``diamond action''), and
hence there is an induced action on
$\ObC^{\omega,m}_{S,Q_m}(\CO)_{\FM_{\rhobar}}$. 
Using Lemma \ref{crucial} and its proof, we easily show: 
\begin{prop}\label{IdentifyCchiAndCm}
We fix a positive integer $m$ and assume the set-up of Assumption~\ref{Aut-BasicGlAss},
with \begin{itemize}
\item[(a)] the set of places $Q_m$ here playing the role of $T$ there,
\item[(b)] $\Lambda$ a discrete valuation ring $\CO$ containing
  $W(\BF)[\zeta_{\ell^m}]$, 
\item[(c)] a set of characters $(\chi_v)$ such that $\chi_v=1$ for
  all $v\in Q_m$, and
\item[(d)] $\ell^m N_n(\ell)$ dividing the order
  of~$k_v^*$ for all $v\in Q_m$.
\end{itemize}
Then $\ObC^{\omega,m}_{S,Q_m}(\CO)_{\FM_{\rhobar}}$ is free
over~$\CO[\Delta_m]$, and for the invariants under $\Delta_m$ one
has
$$\Big(\ObC^{\omega,m}_{S,Q_m}(\CO)_{\FM_{\rhobar}}\Big)^{\Delta_m}=
\ObC^{\omega,0}_{S,Q_m}(\CO)_{\FM_{\rhobar}}.$$
\end{prop}
Note that by Proposition~\ref{Aut-DirSumOfCuspForms} and
Lemma~\ref{stupid}, one has the isomorphism  
$$\ObC^{\omega,m}_{S,Q_m}(\CO)_{\FM_{\rhobar}}\otimes F\cong
\bigoplus_{(\chi_v')}\ObC^{\omega,(\chi_v)}_{S,Q_m}(\CO)_{\FM_{\rhobar}}\otimes
F,$$ 
where the sum is over all $(\chi_v')$ which agree with $(\chi_v)$ for
$v\in S$ and are of order dividing $\ell^m$ for $v\in Q_m$.
\begin{Proof} The second assertion is obvious from the
  definitions. To prove the first, observe that because
  localization at $\FM_{\rhobar}$ commutes with the action of
  $\Delta_m$ and because of Lemma \ref{stupid}, it will suffice to
  show that $\OpbC^{\omega,m}_{S,{Q_m}}(\CO)$ is free
  over~$\CO[\Delta_m]$.

  Because of the central character and the
  conditions at the places $v\in S$, functions in
  $\OpbC^{\omega,m}_{S,{Q_m}}(\CO)$ may be thought of as functions on 
  $$\CD^* Z(K)\backslash\CD^*(\BA)/(\widetilde U_{Q_m}^m
  \prod_{v\in S}\CO^*_{\CD_v}).$$ 
  The $\Delta_m$-orbits of the latter double coset are in bijection with
  the double cosets $\CD^* Z(K)d_j U_{Q_m}^{00}$, $j\in J$, defined in the
  proof of Lemma~\ref{crucial} (with $T={Q_m}$).

Recall that by Proposition~\ref{Aut-DirSumOfCuspForms} one has
\begin{equation}\label{Aut-SumOfSpaces}
\OpbC^{\omega,m}_{S,{Q_m}}(\CO)\otimes_\CO F \cong \bigoplus_{(\chi'_v)}
\OpbC^{\omega,(\chi'_v)}_{S,{Q_m}}(F),
\end{equation}
where the sum is over all character tuples $(\chi'_v)$ such that
$\chi_v=\chi_v'$ for $v\in S$ and $\chi_v'$ is a character of order
dividing $\ell^m$ if $v\in {Q_m}$. Let $u_{j'}$, $j'\in J'$
denote a set of representatives of $\Delta_m$ viewed as a quotient of
$\prod_{v\in {Q_m}}U_0(v)$. The proof of Lemma \ref{crucial} shows that
\begin{enumerate}
\item the values of a cusp form on different double cosets are
  independent,
\item any cusp form is zero on $\CD^* Z(K)d_j U_{Q_m}^{00}$ with
  $j\in J_0$, 
\item and combined with the isomorphism (\ref{Aut-SumOfSpaces}) that functions in
  $\OpbC^{\omega,m}_{S,{Q_m}}(\CO)$ can take arbitrary values on the 
  elements $d_ju_{j'}$, $j\in J\setminus J_0$, $j'\in J'$. \BeweisEnde
\end{enumerate}
\end{Proof}

\smallskip

We now define universal deformation and Hecke rings corresponding to
the above situation in a by now standard fashion regarding the
application of Taylor-Wiles systems:

Choose for each $v\in {Q_m}$ a Teichm\"uller lift
$\lambda_v$ of one of the eigenvalues of $\rhobar(\Frob_v)$. Write
$\Ulambda$ for $(\lambda_v)_{v\in {Q_m}}$. Let
$$\rho_{X,{Q_m}}^{m,\Ulambda}\!:\pi_1(X\setminus (\{w\}\cup{Q_m}))\to \GL_n(R_{X,{Q_m}}^{m,\Ulambda})$$
denote the universal deformation that parameterizes deformations
$\tilde\rho\!:\pi_1(X\setminus (\{w\}\cup{Q_m}))\to\GL_n(R)$, $R\in\CA$, $R$ an
$\CO$-algebra, of $\rhobar$ such that  
\begin{itemize}
\item[(a)] $\det\tilde\rho_w$ factors via $G_w\stackrel{s_w}\longto
  \overline I_w$ and $\tilde\rho_w\otimes(\det\tilde\rho_w)^{-1/n}$ is unramified,
\item[(b)] for all $v\in S$, $\tilde\rho_v(I_v)\to\rhobar(I_v)$ is an isomorphism
\item[(c)] for all $v\in {Q_m}$, $\rho_v$ is of type $\CC_{v,\lambda_v,m}$
  as defined above Lemma~\ref{DependsNotOnV}.
\end{itemize}
Note that the local conditions of $R_{X,{Q_m}}^{m,\Ulambda}$ at $v\in S$ are
described by the local type $\CC_v$ defined in
\cite{boecklekhare},~Prop.~2.15. In particular, at such places one 
has $L_v=H^1_\unr(G_v,\Ad)$. Note also that from the local conditions
at $v\in {Q_m}$, via the action of $I_v$ one obtains a homomorphism
$\CO[\Delta_m]\to R_{X,{Q_m}}^{m,\Ulambda}$.

\smallskip

We define $\omega\!:\GL_1(\BA)\to \CO^*$ as the Hecke-character 
corresponding to $\eta=\det\rho$. Because for $v\in S$ the representation
$\rhobar_v$ is tame and absolutely irreducible and $\rho_v$ is
minimal, we obtain from Lemma~\ref{LocalShapeOfSpecialRho_v} for each $v\in S$ a
character $\chi_v\!:k_{v,n}^*\to W(\BF)^*$, unique up to the Galois
action of $\Gal(k_{v,n}/k_v)$. Finally, for $v\in T$ we let
$\chi_v\!:k_v^*\to \CO^*$ be the trivial character. Then the above data
satisfies Assumption~\ref{Aut-BasicGlAss} for $\Lambda=\CO$. Suppose
a basis of $A_v^n$ is chosen in such a way that $e_{1,v}$ corresponds
to the eigenvalue $\lambda_v$. We now define 
$$\BT_{X,{Q_m}}^{m,\Ulambda}:=(\widetilde\CH_{Q_m}^m(\CO))_{\FM_{\rhobar}}.$$

Let $\tau_1,\ldots,\tau_s$ be a list of the Galois representations
corresponding via Lafforgue's theorem and the global Jacquet-Langlands
correspondence to eigenforms in
$\ObC_{S,{Q_m}}^{m,\omega}(\CO)_{\FM_{\rhobar}}$. By choice of the
maximal ideal and definition of the Hecke action at places in $Q_m$,
using the $\BC$-valued theory one finds that the algebra 
$\BT_{X,{Q_m}}^{m,\Ulambda}\otimes_\CO F$ is semisimple, cf. 
\cite{harristaylor}, (III.2), second paragraph. We therefore
denote by
$$\tau:=\tau_1\oplus\ldots\oplus\tau_s\!:\pi_1(X\setminus
{Q_m})\longto\GL_n(\BT_{X,{Q_m}}^{m,\Ulambda}\otimes_\CO F)$$ 
the corresponding Galois representation.
\begin{prop}
The representation $\tau$ can be written as
the composition of representations
$$\pi_1(X\setminus {Q_m})
\stackrel{\tau_{X,{Q_m}}^{m,\Ulambda}}\longto\GL_n(\BT_{X,{Q_m}}^{m,\Ulambda})
\into\GL_n(\BT_{X,{Q_m}}^{m,\Ulambda}\otimes_\CO F).$$
\end{prop}
\begin{Proof}
Because $\rhobar$ is absolutely irreducible, the image of
$\tau$ lies in the ring of traces. By the
\v{C}ebotarev density theorem, this ring is spanned by the
coefficients of the characteristic polynomials of
$$\tau(\Frob_v), v\notin S\cup\{w\}\cup {Q_m}.$$
Thus by Lafforgue's theorem, it is spanned by the Hecke-eigenvalues of
the corresponding eigenforms. Thus the ring of traces lies in
$\BT_{X,{Q_m}}^{m,\Ulambda}$.
\end{Proof}

From the definition of $\tau_{X,{Q_m}}^{m,\Ulambda}$ it is clear that it is a
representation of the type parameterized by
$\rho_{X,{Q_m}}^{m,\Ulambda}$. By universality there arises a unique morphism
$R_{X,{Q_m}}^{m,\Ulambda}\to \BT_{X,{Q_m}}^{m,\Ulambda}$ such that
$\tau_{X,{Q_m}}^{m,\Ulambda}$ is induced from $\rho_{X,{Q_m}}^{m,\Ulambda}$. 
Both rings are non-zero because of the existence of $\rho$.  

\begin{prop}
The induced morphism $R_{X,{Q_m}}^{m,\Ulambda}\to
\BT_{X,{Q_m}}^{m,\Ulambda}$ is surjective.
\end{prop}
\begin{Proof}
By Nakayama's lemma it suffices to prove the assertion modulo $\FM$. 
It is clear from the Langlands correspondence that the reductions
modulo $\FM$ of the Hecke-operators $T_{v,i}$, $v\notin S\cup
\{w\}\cup T$, and $T'_{w,i}$, $i=1,\ldots,n$, lie in the image. At
places $v\in T$ this again follows from the compatibility of the
global Langlands correspondence with the local one, and 
the explicit decription of the Hecke action at places in $T$,
cf. \cite{harristaylor} (V.1.5).  
Namely, the action of $V_{v,n}\pmod \FM$ is given by the
first eigenvalue $\bar\lambda$ of $\rhobar(\Frob_v)$ and the action of
$V_{v,i}\pmod \FM$, $i=1,\ldots,n-1$, by the elements $b_{v,i}\pmod
\FM$ which are expressions in the elementary symmetric polynomials in
the remaining eigenvalues of $\rhobar(\Frob_v)$. So in this case, too,
the reductions of the Hecke-operators lie in the image of
$R_{X,{Q_m}}^{m,\Ulambda}$.
\end{Proof}

If $Q_m=\emptyset$, we drop it as well as $m$ and $\Ulambda$ from the
notation, and add a superscript zero, i.e., the above morphism becomes
$R_{X}^{0}\to \BT_{X}^{0}$. We have the following central result:
\begin{theorem}\label{TWtheorem}
Suppose that Assumption~\ref{TWsysAss} and the following conditions
hold
\begin{enumerate}
\item For any $\pi_1(X)$-subrepresentation $V$ of
  $\Ad$, there exists a regular semisimple
  $g_V\in\rhobar(\pi_1(\overline X))$ such that $V^{g_V}\neq0$.
\item If $\zeta_\ell\in K$, then $\ad$ has no $1$-dimensional
  subrepresentation,
\item If $\zeta_\ell\in E$, then $H^1(\Gal(E/K(\zeta_\ell),\ad)=0$.
\item $\image(\rhobar)$ contains no normal subgroup of index~$\ell$.
\end{enumerate} 
Then $R_{X}^{0}\to \BT_{X}^{0}$ is an isomorphism. 
\end{theorem}
\begin{Proof}
The proof is based on the use of Taylor-Wiles system in the improved
form due independently to F.~Diamond and K.~Fujiwara, cf.~\cite{diamond} and \cite{fujiwara}. For each $m\in \BN$, use
Lemma~\ref{LemONTWprimes} to choose a set ${Q_m}\subset X$ such that
\begin{itemize}
\item[(a)] $\# Q_m=d:=\dim H^1_{\{L_v\}}(\{w\},\Ad(1))$
\item[(b)] $q_{v}\equiv 1\!\!\pmod{N_n(\ell) \ell^m}$ for all $v\in Q_m$,
\item[(c)] $\rhobar(\Frob_{v})$ has distinct eigenvalues for each $v\in Q_m$, and 
\item[(d)] $H^1_{\{L_{v}^\perp\}}(\{w\}\cup Q_m,\Ad)=0$ where
  $L_v=L_{v,\lambda_{v}}$ for each $v\in Q_m$ and $\lambda_v$ is the
  Teichm\"uller lift of some eigenvalue of $\rhobar(\Frob_{v})$. 
\end{itemize}
Note that $\dim H^1_{\{L_v\}}(\{w\},\Ad)=\dim H^1_{\{L_{v}\}}
(\{w\}\cup Q_m,\Ad)$, and define $\Ulambda_m$ as $(\lambda_v)_{v\in Q_m}$. 

We introduce notation similar to \cite{diamond}, \S 2. Define
\begin{eqnarray*}
{\bf R}&:=&R_{X}^{0}/\FM,\;{\bf T}:=\BT_{X}^{0}/\FM,\\
{\bf R}_m&:=&R_{X,Q_m}^{m,\Ulambda_m}/\FM,\;{\bf T}_m:=\BT_{X,Q_m}^{m,\Ulambda_m}/\FM,\\
{\bf H}&:=&\Hom_{\CO}(\ObC_{S,\emptyset}^{0,\omega}(\CO)_{\FM_\rhobar},\BF),\\
{\bf H}_m&:=&\Hom_{\CO}(\ObC_{S,Q_m}^{m,\omega}(\CO)_{\FM_\rhobar},\BF),\\
{\bf A}_m&:=&\CO[\Delta_m]/\FM\cong\BF[[y_1,\ldots,y_d]]/(y_1,\ldots,y_d)^m 
\end{eqnarray*}
One easily verifies from the preceding work that
\begin{enumerate}
\item Each ${\bf R}_m$ is topologically generated by $d$ elements
  over~$\BF$,
\item ${\bf R}_m/\FM_{{\bf R}_m}^m\cong \BF[[x_1,\ldots,x_d]]/(x_1,\ldots,x_d)^m$,
\item there exists a canonical ${\bf R}_m$-linear surjection $\pi_m\!:{\bf
  H}_m\onto {\bf H}$,
\item under ${\bf R}_m\onto{\bf R}$ the image of $(y_1,\ldots,y_d)$ is
  zero,
\item ${\bf H}_m$ is via ${\bf A}_m\to{\bf R}_m\to{\bf T}_m$ a module
  over ${\bf A}_m$ and ${\bf R}_m$ and the action of ${\bf A}_m$ is
  the same as that which occurs in Proposition~\ref{IdentifyCchiAndCm},
\item ${\bf H}_m$ is free over ${\bf A}_m$ (by Proposition \ref{IdentifyCchiAndCm}).
\end{enumerate}
We now verify the following assertion that is also a crucial part of constructing Taylor-Wiles systems:
\begin{itemize}
\item[(vii)]
The morphism $\pi_m$ induces an isomorphism ${\bf H}_m/(y_1,\ldots,y_d) \cong{\bf H}$ where we consider these as modules over 
$R_{X,Q_m}^{m,\Ulambda_m}$.
\end{itemize}
In view of the second part of Proposition \ref{IdentifyCchiAndCm} and
the above definitions it suffices to show that we have an isomorphism
\begin{equation}\label{TWinclusion}
\ObC_{S,\emptyset}^{0,\omega}(\CO)_{\FM_\rhobar} \simeq
\ObC_{S,Q_m}^{0,\omega}(\CO)_{\FM_\rhobar}
\end{equation}
as modules over $R_{X,Q_m}^{m,\Ulambda_m}$. This follows from  the
arguments in the proofs of \cite{harristaylor}, Prop.~(V.2.3)
and~(V.2.4), and where the isomorphism above is given by the (exact
analog of the) map by $X_{\infty,Q_m}$ of \cite{harristaylor}. We give some
details. As remarked in {\it loc.\ cit.}\ we need to prove the
isomorphy in (\ref{TWinclusion}) only after tensoring with
$\overline{\BQ}_\ell$. We first prove that all forms that contribute
to the right hand side are old at places in $Q_m$. 

For this
let $\tilde f$ be a form in the latter space and let $f$ be the
corresponding cusp form for $\GL_n$ under the Jacquet-Langlands
correspondence (that is applicable here after the work of Badulescu as
$\tilde f$ is supercuspidal at a sufficiently large number of
places). Let $v$ be in $Q_m$. As the $\ell$-adic representation
$\rho_f$ corresponding to $f$ reduces to $\rhobar$ residually,
$\rho_f$ restricted to a decomposition group $G_v$ at $v$ is a lift of
$\rhobar_v$ to $\CO$. By Lemma~\ref{LemmaOnTWPrimes} such lifts are
diagonalizable and therefore finite on inertia. Hence by the
compatibility of the global Langlands correspondence of
\cite{lafforgue} with the local Langlands correspondence and because
of triviality of the ``nebentype'', $\Pi_v(f)$ is an unramified principal series representation.

We are now in the situation described on page~\pageref{HTquote} and
quoted from \cite{harristaylor}. The action of $V_{v,n}$ on
$\Pi_v(f)^{U_0(v)}$ is diagonalizable and there is exactly one
eigenvalues whose reduction modulo $\FM$ is $\lambda_v\pmod \FM$, and
it has multiplicity one. Therefore after localization at
$\FM_{\rhobar}$ at most a one-dimensional subspace of
$\Pi_v(f)^{U_0(v)}$ remains. That the remaining space is indeed
$1$-dimensional follows from the fact that (\ref{TWinclusion}) is
injective. We have thus completed the proof of~(vii).

\smallskip

It now follows from \cite{diamond}, Thm~2.1, that ${\bf R}$  is a
complete intersection of dimension zero and ${\bf H}$ is free over  
${\bf R}$. Since the action of ${\bf R}$ on ${\bf H}$ factors via
${\bf T}$, it follows in particular that ${\bf R}\to{\bf T}$ must be
injective and hence an isomorphism. 

Let us now come back to the original question. Because
${\BT_{X,Q_m}^{m,\Ulambda}}$ is $\CO$-torsion free and finitely
generated, the surjection
$R_{X,Q_m}^{m,\Ulambda}\onto\BT_{X,Q_m}^{m,\Ulambda}$ splits as
a map of $\CO$-modules. By the above its reduction modulo $\FM$ is an
isomorphism. But then the morphism itself must be bijective. 
\end{Proof}

\subsection{Lowering the level \`a la Skinner and Wiles}\label{LevelLow}

We have the following `lowering the level' result 
as in the work of Skinner and Wiles in \cite{skinnerwiles}.
\begin{theorem}\label{GalRepLoweringLevel}
Suppose $\rho\!:\pi_1(X\setminus T)\to\GL_n(\CO)$ and
$\rhobar:=\rho\pmod {\FM}$ satisfy the
following conditions:
\begin{itemize}\item[(a)] $\CO$ contains $\zeta_\ell$,
\item[(b)]  $\# S\geq 2$, and for all $v\in S$ the representation
  $\rhobar_v$ is absolutely irreducible and tamely ramified, and
  $\rho_v$ is minimal,
\item[(c)] for all $v\in T$, $\rho$ is of type-$1$ and $\ell N_n(\ell)$ divides the order
  of~$k_v^*$.
\end{itemize}

Then there exists a
representation $\rho'\!:\pi_1(X\setminus T)\to \GL_n(\CO)$ such that
\begin{enumerate}
\item the residual representations of $\rho$ and $\rho'$ agree,
\item $\rho_v'$ is minimal for $v\in S$ and
\item $\rho'(I_v)$ is finite for $v\in T$.
\end{enumerate} 
\end{theorem}

\begin{remark} {\em This theorem is referred to as a level lowering result as 
from it one deduces that there is a solvable base change $Y \rightarrow X$, that 
one can make disjoint from any given covering of $X$, such that $\rho'|_{\pi_1(Y)}$
has conductor the conductor of $\rhobar|_{\pi_1(Y)}$. We use this in the section that is coming up!}
\end{remark}

\begin{Proofof}{Theorem~\ref{GalRepLoweringLevel}} We use Lafforgue's
  theorem to convert the above into an assertion  about cusp
  eigenforms that we have proved in Lemma~\ref{AutoLoweringLevel}: Via
  Lafforgue's theorem, which is compatible with the local Langlands
  correspondence, $\rho$ corresponds to a cuspidal Hecke eigenform in
  $\bC^{\omega,(\chi_v)}_{S,T}(\CO)_{\FM_{\rhobar}}$. Badulescu's
  Theorem~\ref{GlobalThmBadu} in the form of
  Corollary~\ref{SecondCorOnBadu} yields a corresponding Hecke
  eigenform in the space
  $\ObC^{\omega,(\chi_v)}_{S,T}(\CO)_{\FM_{\rhobar}}$ (here we do not need to consider Hecke action at places in $T$).

Let now $(\chi'_v)$ be a set of characters which agrees with
$(\chi_v)$ for $v\in S$ and is of exact order $\ell$ at all $v\in
T$. By Lemma~\ref{AutoLoweringLevel} and again
Corollary~\ref{SecondCorOnBadu} 
we therefore find a non-zero cuspidal Hecke eigenform $f'$ in
$\bC^{\omega,(\chi'_v)}_{S,T}(\CO)_{\FM_{\rhobar}}$. 
Let $\rho'$ be the corresponding Galois representations, which exists
by Lafforgue's theorem. Assertion (i) is clear from the definition of
$\FM_{\rhobar}$. At places in $T$, it follows from the non-triviality
of $\chi'_v$, that $\Pi_v(f')$ is ramified principal series, 
i.e, potentially unramified, because at such the character $\chi'_v$
is non-trivial. This shows (iii). Finally (ii) follows
from the definition of $\ObC^{\omega,(\chi'_v)}_{S,T}(F)$
together with Propositions~\ref{LocalBadulescu}
and~\ref{LocalLanglForChi}.
\end{Proofof}

\section{Proof of main theorems}\label{ThmBProof}

We can finally give the proof of our central theorem.
\begin{Proofof}{Theorem~\ref{OnDeJConj}}
We argue by contradiction, and, in view of Theorem~\ref{ConjEqConj},
assume that there is a representation
$\rho'\!:\pi_1(X)\to\GL_n(\BF[[x]])$ with $\rhobar =\rho'\!\!\pmod x$
and $\rho'(\pi_1(\overline X))$ infinite. Because $\ell$ does not
divide $n$, taking $n$-th roots is an isomorphism on the $1$-units in
$\BF[[x]]^*$. Therefore the $n$-th root of the $1$-unit part of
$\det\rho'$ is a character, say $\tilde\eta$, of $\pi_1(X)$. It has the
property $\det(\tilde\eta\cdot\rho')=\det\rhobar$. Since the image of
$\pi_1(\overline X)$ under any character of $\pi_1(X)$ is finite, we
will from now on assume $\det\rho'=\det\rhobar$.

Let $v$ be a place as in \ref{CondFive}. We claim that $\rho'$ is
minimal at $v$: Note first that ramification at $v\in S$ is tame since
$\ell\neq p$ and $\rhobar_v$ is tame. Now the tame inertia is pro-cyclic and
by Lemma~\ref{LocalShapeOfSpecialRho_v} the representation
$\rhobar_{|I_v}$ is diagonal with distinct characters. Therefore
$\rho'_{|I_v}$ is diagonal with distinct characters and of order
dividing $q_v^n-1$. Because $(1+x\BF[[x]],\cdot)$ is torsion free, we
must have $\rho'(I_v)\cong\rhobar(I_v)$.

\smallskip

We now consider the ring $R:=\kernel(\CO\oplus\BF[[x]]\onto\BF)$. It
lies in $\CA$ and affords a representation
$\rho'':=\rho\oplus\rho'\!:\pi_1(X)\to\GL_n(R)$ with determinant
$\eta$. By \cite{dejong}, Lemma~2.12, the ramification of $\rho'$ at
all places in $S$ is finite. Let $v$ be a place at 
which $\rhobar_v$ is tamely ramified and absolutely irreducible and
such that $\rho(I_v)\cong\rhobar(I_v)$. Set $T_s:=\{v\}$ and $T_i:=T$,
the set of places at which $\rho$ is of type-$1$. Take $m\in\BN$ such
that $\ell N_n(\ell)$ divides $\# k_{w,m}^*$ for all $w\in T_i$.
Then Corollary~\ref{EnlargeBaseFieldCor} applied to $\rho''$ provides
us with a finite Galois covering $Y\to X$ such that
\begin{enumerate}
\item $Y$ is geometrically connected over $k$. 
\item $\rho''(\pi_1(Y))=\rho''(\pi_1(X))$, $\rhobar(\pi_1(\overline
  Y))=\rhobar(\pi_1(\overline X))$,
\item $\rho'_{|\pi_1(Y)}$ is ramified precisely at the places above
  $v$, 
\item $\rhobar_{v'}$ is tame and absolutely irreducible and $\rho''_v$
  is minimal at places $v'$ above $v$, 
\item $\rho$ is unramified at places not above $T_i\cup\{v\}$, and
\item $\rho_{|G_{w'}}$  is of type-$1$ and $\ell\cdot
  N_n(\ell)|\#k_{w'}^*$ for all $w'$ above a place $w\in T_i$.
\end{enumerate}
By passing to a second cover, if necessary, we may furthermore assume
that there are at least two places in $Y$ above~$v$.

Since $\rho_{\pi_1(Y)}$ satisfies all the conditions originally
imposed on $\rho$, we may therefore rename $Y$ to $X$, assume that
$\rho'$ is unramified outside $S$ and has determinant $\det\rhobar$,
and that $\rho$ satisfies the following conditions: 
\begin{enumerate}
\item $\rhobar:=\rho\!\!\pmod{\FM}$ is absolutely irreducible,
\item $\rhobar(\pi_1(\overline X))$ contains a regular semisimple element,
\item $\eta:=\det\rho$ is of finite order,
\item \label{CondFour}at places $v$ at which $\rho$ ramifies and
  $\rhobar$ does not, $\rho_v$ is of type-$1$ and $\ell N_n(\ell)$ divides
  $\# k_{v}^*$, 
\item at places $v$ at which $\rhobar$ ramifies, $\rhobar_v$ is tame and
  absolutely irreducible and $\rho''_v$ is minimal, and
\item $\rhobar$ ramifies at least at two places.
\end{enumerate}

Condition~\ref{CondFour} allows us to apply
Theorem~\ref{GalRepLoweringLevel} on level lowering. Thereby we may
replace (iv) by
\begin{itemize}
\item[{\ref{CondFour}'}] $\rho(I_v)$ is finite at all
  places where $\rho$ ramifies and $\rhobar$ does not.  
\end{itemize}

By yet another application of Corollary~\ref{EnlargeBaseFieldCor} to
$\rho''$ with $T_i=\emptyset$ and $T_s$ the places where $\rhobar$
ramifies, the latter condition may be replaced by 
\begin{itemize}
\item[{\ref{CondFour}''}] $\rho''$ is ramified only at places at which
  $\rhobar$ is.
\end{itemize}

By enlarging $X$, we now assume that $S$ is the set of places at which
$\rhobar$ is ramified. Recall that $R_X^{0}$ is universal for
deformations of $\rhobar$ to $\CO$-algebras in $\CA$ which are
unramified outside $S\cup\{w\}$, minimal at $S$ and unramified at $w$
after a twist by a character of order $\ell$ that factors via
$s_w$. Because $\rho'$ is minimal at $S$ and unramified outside $S$,
there is a unique morphism $\phi\!:R_X^{0}\to \BF[[x]]$ which
induces $\rho'$. The ring $\BF[[x]]$ is of characteristic $\ell$ and
so $\phi$ factors via~$R_X^{0}/(\ell)$.

On the other hand by Theorem~\ref{TWtheorem}, which used the technique
of Taylor-Wiles systems, the ring $R_X^{0}$ is finite over
$\BZ_\ell$. Therefore $R_X^{0}/(\ell)$ is finite, and this
contradicts our assumption that $\rho'$ has infinite image. 
\end{Proofof}

When combined with Lemma~\ref{SL_nAdmitsTWclasses} and
\cite{boecklekhare},~Prop.~2.3, the following result implies
Theorem~\ref{SpecialDeJConj}. (Note that (i) implies $\ell\notdiv n$.)
\begin{theorem}\label{OnDeJConj2}
Suppose $\rhobar\!:\pi_1(X)\to \GL_n(\BF)$ satisfies
\begin{enumerate}
\item $\ad$ is absolutely irreducible over $\BF_\ell[\image(\rhobar)]$,
\item if $\zeta_\ell\in E$, then $H^1(\Gal(E/K(\zeta_\ell)),\ad)=0$
  and $\ad$ is absolutely irreducible over $\BF[\rhobar(\pi_1(Z))]$,
  where $Z\to X$ corresponds to $K(\zeta_\ell)/K$,
\item $\rhobar$ admits $R$-places, 
\item $\rhobar(\pi_1(\overline X))$ contains a regular semisimple element, and
\item there exists a place $v\in S$ such that $\rhobar_{v}$ is tame
  and absolutely irreducible.
\end{enumerate}
Then $R^\eta_{\rhobar}$ is finite over~$\BZ_\ell$, where $\eta$ is the
Teichm\"uller lift of $\det\rhobar$.
\end{theorem}
\begin{Proof}
To prove the assertion on $R_X^\eta$ we may, as in the preceeding
proof, pass from $X$ to a finite Galois cover provided that we preserve
all our original hypothesis. Therefore by
Corollary~\ref{EnlargeBaseFieldCor} we may replace condition~(v) by  
\begin{itemize}
\item[(v)'] at all places $v$ at which $\rhobar$ is ramified,
  $\rhobar_v$ is absolutely irreducible and tame.
\end{itemize}

Using (i)--(iii) and (v') we obtain from
\cite{boecklekhare},~Thm.~2.1, a finite set $T\subset X$ and a
representation $\rho\!:\pi_1(X\setminus T)\to\GL_n(W(\BF))$ such that
$\rho$ is of type-$1$ at places in $T$ and minimal above places
in~$S$.

Suppose now that $\image(\rhobar)$ contains a normal subgroup of index
$\ell$, and let $\pi\!:Y\to X$ be the corresponding Galois cover of
degree $\ell$. Because $\ell\notdiv n$, the modular representation
theory of finite groups shows that $\ad$ is still absolutely
irreducible over $\rhobar(\pi_1(Y))$. Also (iv) and (v) still hold for
$\rhobar|_{\pi_1(Y)}$. Moreover the splitting field over $k(Y)$ of
$\rhobar$ is still~$E$. We claim that (ii) still holds for
$\rhobar_{|\pi_1(Z)}$.

So suppose $\zeta_\ell\in E$. By the reasoning given above, $\ad$ will still be absolutely
irreducible over $\BF[\pi_1(Z_Y)]$ for the pullback $Z_Y\to Y$ of
$Z\to X$ along $Y\to X$. Let $K'/K(\zeta_\ell)$ be the field extension
corresponding to $Z_Y\to Y$. Inflation-restriction yields
\begin{eqnarray*}\lefteqn{0\to H^1(\Gal(K'/K(\zeta_\ell)),\ad^{\Gal(E/K')})\to
H^1(\Gal(E/K(\zeta_\ell),\ad)}\\
&\!\!\!\!\to& \!\!\!H^1(\Gal(E/K',\ad)^{\Gal(K'/K(\zeta_\ell))}\to
H^2(\Gal(K'/K(\zeta_\ell)),\ad^{\Gal(E/K')}).
\end{eqnarray*}
The outer terms are zero because $\ad^{\Gal(E/K')}=0$.
The second term is zero by assumption. Now any $\ell$-group acting on
a finite-dimensional non-trivial $\BF_\ell$ vector space has a
non-trivial set of invariants. Since
$\Gal(K'/K(\zeta_\ell))\cong\BZ/(\ell)$,  this implies
$H^1(\Gal(E/K',\ad)=0$. Thus (ii) holds over
$Y$ instead of $X$. 

By induction, we may therefore pass to an
extension $Y'$ of $X$ over which (i), (ii), (iv) and (v) hold, and
such that in addition $\rhobar(\pi_1(Y'))$ has no normal subgroup of
index $\ell$. Over $Y'$ we can now apply Theorem~\ref{OnDeJConj}, and
the result follows. (This uses again the formulation of
Conjecture~\ref{DeJConj}, which makes it obvious that de Jong's
conjecture holds for $\rhobar$, if it holds for $\rhobar_{|\pi_1(Y')}$.)
\end{Proof}

\section{Appendix}

In the paper \cite{taylorwiles} a method was developed to show that Hecke rings are complete 
intersections in ``minimal cases''. This was done via constructing what are now called Taylor-Wiles 
(TW-) systems. Since then there have been many works that have made technical simplifications to the
method of \cite{taylorwiles}. For instance in \cite{diamond} and
\cite{skinnerwiles} notable simplifications are achieved. 
There is a trick in the main body of the paper that produces another
(rather minor admittedly, but perhaps useful!) simplification which pertains to
the avoidance of imposing an auxiliary level structure to handle ``torsion problems''.
In this appendix we explain this trick in the original context of modular curves of \cite{taylorwiles}.
We use notation of \cite{diamond} to indicate what the problem is and how we handle it. The main point is that when 
proving freeness of certain cohomology groups it is enough for the purposes of TW systems to prove this
over certain group algebras whose group of characters ``kill torsion''.

A key step in TW systems is to prove, for certain finite set of primes $Q=\{q_1,\cdots,q_r\}$
and any positive integer $N$ prime to the primes in $Q$, results towards the freeness of 
the cohomology group $H^1(X_{N,Q},{\cal O})_{\sf m}$ as 
a module over ${\cal O}[\Delta_Q]$ (under the natural action) with $\Delta_Q$ the Sylow $\ell$-subgroup 
(which we may also view as the maximal $\ell$-quotient) of $\Pi_{q \in Q} ({\bf Z}/q{\bf Z})^*$,
where $\sf m$ is a mod $\ell$ maximal ideal of a certain Hecke algebra which satisfies a certain set of conditions, and where $\cal O$ 
is  a  finite flat extension of ${\bf Z}_{\ell}$.
Here $X_{N,Q}$ is the modular curve corresponding to the subgroup $\Gamma_0(N) \cap \Gamma_1(q_1 \cdots q_r)$. Further
the quotient of $H^1(X_{N,Q},{\cal O})_{\sf m}$ by the augmentation ideal of ${\bf Z}_{\ell}[\Delta_Q]$ is isomorphic
to $H^1(X_0(N),{\cal O})_{\sf m}$ if the primes $q_i$ in addition 
satisfy the hypothesis that the mod $\ell$ representation $\rhobar$ corresponding to $\sf m$ is such that it is unramified at $q_i$ and the ratio of the eigenvalues of $\rhobar({\rm Frob}_{q_i})$ are not $q_i^{\pm 1}$.

One of the steps in proving the freeness is to impose a {\it usefully harmless} auxiliary level structure, 
{\it useful} in that it gets rid of torsion in the congruence subgroup, 
and {\it harmless} in that it does not introduce new forms that 
give rise to the same residual representation that is fixed upon, to avoid problems arising from torsion of 
$\Gamma_0(N)$ (or better still, the smaller torsion of $\Gamma_0(Nq_1\cdots q_r)$). 
There are situations where such an auxiliary level structure may not be found, and thus we indicate in this setting 
a more elementary trick  developed in the main body of the paper (which was in the case of anisotropic groups) to proving such freeness results without any need
of imposing auxiliary level structure that are enough in applications to build TW systems using the first cohomology of modular curves, 
sketching very briefly the proofs
(which are completely standard). (As choosing an auxiliary level structure uses somewhat delicate information about
the residual representation and the compatibility of the local and global Langlands correspondence at ramified places, we regard our modification as more elementary, and yet another instance of the flexibility that inheres in the Taylor-Wiles method.)
The cohomology groups we consider below can either be Betti or \'etale.

We observe the following proposition which directly follows from the proof of \cite{carayol}~Lemme~1:

\begin{propA}
  Let $e$ be any integer that kills the torsion of $\Gamma_0(N)$. Let $\Delta_e^Q$ be 
  the subgroup  of ${\rm Hom}(\Delta_Q,{{\bf Q}_{\ell}/{\bf Z}_{\ell}})$ that consists of $e$th powers.
   Consider the twisted sheaf
  ${\cal O}(\chi)$ for any character $\chi \in \Delta_e^Q$ on $X_0(Nq_1 \cdots q_r)$ and assume ${\cal O}$ to be large enough to contain all values of $\chi$. Then if $k$ is the residue field of
  ${\cal O}$, the reduction ${\cal O}(\chi) \otimes k$ is isomorphic to the constant sheaf $k$ on $X_0(Nq_1 \cdots q_r)$.
\end{propA}

As a standard consequence one has: 
\begin{corA} Let $\Delta_Q^e$ be the quotient of $\Delta_Q$
  that is dual to the subgroup $\Delta_e^Q$ of ${\rm Hom}(\Delta_Q,{{\bf Q}_{\ell}/{\bf Z}_{\ell}})$.
 Let $X_{N,Q,e}$ be the modular curve that corresponds to the congruence subgroup that is the 
kernel of the natural map $\Gamma_0(Nq_1 \cdots q_r) \rightarrow \Delta_Q^e$. Let $\sf m$ be 
a maximal ideal of a certain Hecke algebra (as in \cite{diamond}: we drop operators $T_r,U_r$ for $r$ 
not coprime to $Nq_1\cdots q_n$) acting on $H^1(X_{N,Q,e},{\cal O})$
  such that the corresponding (semisimple by definition) residual representation $\rhobar$ of $G_{\bf Q}$ is not abelian. Then
  $H^1(X_{N,Q,e},{\cal O})_{\sf m}$ has a natural action of $\Delta_Q^e$ and is a free ${\cal O}[\Delta_Q^e]$-module for any ${\cal O}$
  that is finite flat over ${\bf Z}_{\ell}$.
\end{corA} 

\begin{Proofof}{the proposition}
For conciseness of notation we denote by 
$Y$ and $X$ the curves $X_{N,Q,e}$ and $X_0(Nq_1 \cdots q_r)$, and we have the natural map $r:Y \rightarrow X$ that is the quotient
by $\Delta_Q^e$.
The sheaf ${\cal O}(\chi)$ is described as $\Delta_Q^e \backslash [Y \times {\cal O}]$ with $\Delta_Q^e$ acting
on the constant sheaf ${\cal O}$ by $\chi$. The stalk at a point $x \in X$, after choosing a point $y$ in $r^{-1}(x)$,
can be identified with the subset of 
the stalk at $y$ of the constant sheaf, ${\cal O}_y$, on which the stabilizer of $y$ in $\Delta_Q^e$
acts by $\chi$ (thus it is either ${\cal O}_y$ or 0). 
From this description the proposition follows using our assumption on $\chi$.\end{Proofof}

\begin{Proofof}{the corollary} We first note that as $\rhobar$ is irreducible, the \'etale $H^0$ and $H^2$ 
of modular curves with coefficients in the twisted sheaves above
do not have the maximal ideal $\sf m$ in their support (see Section 3 of \cite{carayol}). Thus from the proposition above, and the long exact sequence 
of cohomology, it follows  that for 
 for any character $\chi \in \Delta_e^Q$ on $X_0(Nq_1 \cdots q_r)$ we have a (Hecke equivariant) isomorphism
$H^1(X_0(Nq_1 \cdots q_r),{\cal O}(\chi))_{\sf m} \otimes k  \simeq H^1(X_0(Nq_1 \cdots q_r,k)_{\sf m}$.
Then by a standard argument (see proof of Proposition 5.6.1 in \cite{cdt}) the corollary follows.
\end{Proofof}

Let us further assume that for each $n\in\BN$ we have sets of primes
$Q_n$ of constant cardinality $r$ such that for $q \in Q_n$, $q$ is
prime to $N$,  $q$ is 1 mod $\ell^n$, and  $\rhobar$ is unramified at
$q$ with the ratio of the eigenvalues of $\rho({\rm Frob}_{q})$ not
$q^{\pm  1}$. Then again $H^1(X_{N,Q,e},{\cal O})_{\sf m}$ is a free
${\cal O}[\Delta_{Q_n}^e]$-module whose quotient by $\Delta_{Q_n}^e$ 
is isomorphic to $H^1(X_0(N),{\cal O})_{\sf m}$. The group
$\Delta_{Q_n}^e$ surjects onto $\BZ/(\ell^{n-e})^r$, and thus grows
systematically with $n$. This is enough to construct TW systems
as in Section 3.1  of \cite{diamond} avoiding imposition of auxiliary
level structures.

The trick of this appendix can also be used to avoid imposition of
auxiliary level structure needed to bypass torsion problems in the
level-lowering method of \cite{skinnerwiles}, as done in the main body
of this paper. Here the further remark, in addition to the observation above, is 
that when base changing to make orders of certain residue fields congruent
to 1 modulo high powers of $\ell$, one also requires certain other
places, chosen in advance, to split completely, so that the
$\ell$-part of the torsion can not grow much under base change.

\noindent {\it Addresses of the authors}: 

\noindent GB: Institut f\"ur experimentelle Mathematik, Universit\"at
Duisburg-Essen, Standort Essen, Ellernstrasse 29, 45326 Essen, Germany. \\
e-mail address: {\tt boeckle@exp-math.uni-essen.de}

\noindent CK: Dept of Math, University of Utah, 155 S 1400 E,
Salt Lake City, UT 84112.\\ 
e-mail address: {\tt shekhar@math.utah.edu}
 
\noindent School of Mathematics, 
TIFR, Homi Bhabha Road, Mumbai 400 005, INDIA.\\ 
e-mail addresses: {\tt shekhar@math.tifr.res.in}
 
\end{document}